\documentclass[matprg]{mcsreport}

\usepackage{subeqn}
\usepackage{extra}
\usepackage{rotating}
\usepackage{epsfig}

\def\labtag#1{\label{#1}}

\newcommand{\epstol}{\epsilon_{\rm tol}}

\begin{document}
\author{Jeff Linderoth \and Stephen Wright}

\title{Decomposition Algorithms for Stochastic
  Programming on a Computational Grid}

\titlerunning{Stochastic Programming on a Computational Grid}

\institute{Jeff Linderoth\at
Axioma Inc., 501-F Johnson Ferry Road, Suite 450,
Marietta, GA 30068;
{\tt jlinderoth@axiomainc.com}
\and 
Stephen Wright\at 
Mathematics and Computer Science Division,
Argonne National Laboratory, 9700 South Cass Avenue, 
Argonne, IL 60439; {\tt wright@mcs.anl.gov}}
\date{\today}
\reportnumber{P875--0401, April, 2001}
\maketitle

\begin{abstract}
  We describe algorithms for two-stage stochastic linear programming
  with recourse and their implementation on a grid computing platform.
  In particular, we examine serial and asynchronous versions of the
  L-shaped method and a trust-region method. The parallel platform of
  choice is the dynamic, heterogeneous, opportunistic platform
  provided by the Condor system. The algorithms are of master-worker
  type (with the workers being used to solve second-stage problems),
  and the MW runtime support library (which supports master-worker
  computations) is key to the implementation.  Computational results
  are presented on large sample average approximations of problems
  from the literature.
\end{abstract}

\section{Introduction} \labtag{introduction}

Consider the following stochastic optimization problem:
\beq \labtag{gen.sp}
\min_{x \in S} \, F(x) \defeq \sum_{i=1}^N p_i f(x,\omega_i),
\eeq
where $S \in \R^n$ is a constraint set, $\Omega = \{ \omega_1,
\omega_2, \dots, \omega_N \}$ is the set of outcomes (consisting of
$N$ distinct scenarios), and $p_i$ is the probability associated with
each scenario. Problems of the form \eqnok{gen.sp} can arise directly
(in many applications, the number of scenarios is naturally finite),
or as discretizations of problems over continuous probability spaces,
obtained by approximation or sampling. In this paper, we discuss the
{\em two-stage stochastic linear programming problem with fixed
  resource}, which is a special case of \eqnok{gen.sp} defined as follows:
\begin{subequations} \labtag{2stage.lp}
\beqa 
\labtag{2stage.lp.obj}
& \min \, c^T x + \sum_{i=1}^N p_i q(\omega_i)^T y(\omega_i), \sgap
\mbox{subject to} \\
\labtag{2stage.lp.x}
& Ax=b, \;\; x \ge 0, \\
\labtag{2stage.lp.y}
& W y(\omega_i) = h(\omega_i) - T(\omega_i) x, \;\; 
y(\omega_i) \ge 0, 
\sgap i=1,2,\dots,N.
\eeqa
\end{subequations}
The unknowns in this formulation are $x$ and $y(\omega_1),
y(\omega_2), \dots, y(\omega_N)$, where $x$ contains the ``first-stage
variables'' and each $y(\omega_i)$ contains the ``second-stage
variables'' associated with the $i$th scenario. The $i$th scenario is
characterized by the probability $p_i$ and the data objects
$(q(\omega_i), T(\omega_i), h(\omega_i))$.

The formulation \eqnok{2stage.lp} is sometimes known as the
``deterministic equivalent'' because it lists the unknowns for all
scenarios explicitly and poses the problem as a (potentially very
large) structured linear program. An alternative formulation is
obtained by recognizing that each term in the second-stage summation
in \eqnok{2stage.lp.obj} is  a piecewise linear convex function
of $x$. Defining the $i$th second-stage problem as a linear program (LP)
parametrized by the first-stage variables $x$, that is, 
\begin{subequations}
\labtag{second-stage-lp}
\beqa
\labtag{second-stage-lp.1}
& \cQ_i(x) \defeq \min_{y(\omega_i)} \,  q(\omega_i)^T y(\omega_i) \;\; 
\mbox{subject to} \\
\labtag{second-stage-lp.2}
&  W y(\omega_i) = h(\omega_i) - T(\omega_i) x, 
\;\;  y(\omega_i) \ge 0,
\eeqa
\end{subequations}
and defining the objective in \eqnok{2stage.lp.obj} as
\beq \labtag{def.Q}
\cQ(x) \defeq c^Tx + \sum_{i=1}^N p_i \cQ_i(x),
\eeq
we can restate \eqnok{2stage.lp} as 
\beq \labtag{2stage.pl}
\min_x \, \cQ(x), \;\; \mbox{subject to} \; Ax=b, \; x \ge 0.
\eeq

We note several features about the problem \eqnok{2stage.pl}.  First, it
is clear from \eqnok{def.Q} and \eqnok{second-stage-lp} that $\cQ(x)$
can be evaluated for a given $x$ by solving the $N$ linear programs
\eqnok{second-stage-lp} separately. Second, we can derive subgradient
information for $\cQ_i(x)$ by considering dual solutions of
\eqnok{second-stage-lp}. If we fix $x=\hat{x}$ in
\eqnok{second-stage-lp}, the primal solution $y(\omega_i)$ and dual
solution $\pi(\omega_i)$ satisfy the following optimality conditions:
\beqas
q(\omega_i) - W^T \pi(\omega_i) \ge 0 & \perp &  y(\omega_i) \ge 0, \\
W y(\omega_i) & = & h(\omega_i) - T(\omega_i) \hat{x}.
\eeqas
From these two conditions we obtain that
\beq \labtag{theta.2}
\cQ_i(\hat{x}) = q(\omega_i)^T y(\omega_i) = 
\pi(\omega_i)^T W y(\omega_i) =
\pi(\omega_i)^T [ h(\omega_i) - T(\omega_i) \hat{x} ].
\eeq
Moreover, since $\cQ_i$ is piecewise linear and convex, we have for
any $x$ that
\beq \labtag{subg.property}
\cQ_i(x) - \cQ_i(\hat{x}) \ge 
\pi(\omega_i)^T [ -T(\omega_i) x  + T(\omega_i) \hat{x} ] =
\left( - T(\omega_i)^T \pi(\omega_i) \right)^T (x-\hat{x}),
\eeq
which implies that
\beq \labtag{subg.Qi}
-T(\omega_i)^T \pi(\omega_i) \in \partial \cQ_i(\hat{x}),
\eeq
where $\partial \cQ_i(\hat{x})$ denotes the subgradient of $\cQ_i$ at
$\hat{x}$. By Rockafellar~\cite[Theorem~23.8]{Roc70}, using
polyhedrality of each $\cQ_i$, we have from \eqnok{def.Q} that
\beq \labtag{subg.Q}
\partial \cQ(\hat{x}) = c + \sum_{i=1}^N p_i \partial \cQ_i(\hat{x}),
\eeq
for every $\hat{x}$ that lies in the domain of each $\cQ_i$,
$i=1,2,\dots,N$.

Let $\cS$ denote the solution set for \eqnok{2stage.pl}; we assume for
most of the paper that $\cS$ is nonempty. Since \eqnok{2stage.pl} is a
convex program, $\cS$ is closed and convex, and the projection
operator $P(\cdot)$ onto $\cS$ is well defined.  Because the objective
function in \eqnok{2stage.pl} is piecewise linear and the constraints
are linear, the problem has a {\em weak sharp minimum} (Burke and
Ferris~\cite{BurF93}); that is, there exists $\hat{\epsilon}>0$ such
that
\beq \labtag{weak.sharp}
\cQ(x) - \cQ^* \ge \hat{\epsilon} \| x- P(x) \|_{\infty}, 
\;\; \mbox{for all $x$ with  $Ax=b$, $x \ge 0$,}
\eeq
where $\cQ^*$ is the optimal value of the objective.

The subgradient information can be used by algorithms in different
ways. Successive estimates of the optimal $x$ can be obtained by
minimizing over a convex underestimate of $\cQ(x)$ constructed from
subgradients obtained at earlier iterations, 
as in the L-shaped method described in
Section~\ref{sec:lshaped}. This method can be stabilized by the use of
a quadratic regularization term (Ruszczy{\'n}ski~\cite{Rus86},
Kiwiel~\cite{Kiw90}) or by the explicit use of a trust region, as in
the $\ell_{\infty}$ trust-region approach described in
Section~\ref{sec:tr}.  Alternatively, when an upper bound on the
optimal value $\cQ^*$ is available, one can derive each new iterate
from an approximate analytic center of an approximate epigraph. The
latter approach has been explored by Bahn et al.~\cite{BahDGV95} and
applied to a large stochastic programming problem by Frangi{\`e}re,
Gondzio, and Vial~\cite{FraGV00}.

Because evaluation of $\cQ_i(x)$ and elements of its subdifferential can be
carried out independently for each $i=1,2,\dots,N$, and because such
evaluations usually constitute the bulk of the computational workload,
implementation on parallel computers is possible.  We can partition
second-stage scenarios $i=1,2,\dots,N$ into ``chunks'' and define a
computational task to be the solution of all the LPs
\eqnok{second-stage-lp} in a single chunk. Each such task could be
assigned to an available worker processor. Relationships between the
solutions of \eqnok{second-stage-lp} for different scenarios can be
exploited within each chunk (see Birge and
Louveaux~\cite[Section~5.4]{BirL97}).  The number of second-stage LPs
in each chunk should be chosen to ensure that the computation does
not become communication bound. That is, each chunk should be large
enough that its processing time significantly exceeds the time
required to send the data to the worker processor and to return the
results.

In this paper, we describe implementations of decomposition algorithms
for stochastic programming on a dynamic, heterogeneous computational
grid made up of workstations, PCs (from clusters), and supercomputer
nodes.  Specifically, we use the environment provided by the Condor
system~\cite{condor}.  We also discuss the MW runtime library (Goux et
al.~\cite{GouLY00,GouKLY00}), a software layer that significantly
simplifies the process of implementing parallel algorithms in Condor.

For the dimensions of problems and parallel platforms considered in
this paper, evaluation of the functions $\cQ_i(x)$ and their
subgradients at a single $x$ often is insufficient to make
effective use of the available processors. Moreover, ``synchronous''
algorithms---those that depend for efficiency on all tasks completing
in a timely fashion---run the risk of poor performance in an
environment such as ours, in which failure or suspension of worker
processors while they are processing a task is not an infrequent
event.  We are led therefore to ``asynchronous'' approaches that
consider different points $x$ simultaneously.  Asynchronous variants
of the L-shaped and $\ell_{\infty}$ trust-region methods are described
in Sections~\ref{sec:lshaped:async} and \ref{sec:atr}, respectively.

Other parallel algorithms for stochastic programming have been devised
by Birge et al.~\cite{BirDHS98}, Birge and Qi~\cite{BirQ88}, and
Frangi{\`e}re, Gondzio, and Vial~\cite{FraGV00}.  In \cite{BirDHS98}, the
focus is on multistage problems in which the scenario tree is
decomposed into subtrees, which are processed independently and in
parallel on worker processors. Dual solutions from each subtree are
used to construct a model of the first-stage objective (using an
L-shaped approach like that described in Section~\ref{sec:lshaped}),
which is periodically solved by a master process to obtain a new
candidate first-stage solution $x$.  Parallelization of the linear
algebra operations in interior-point algorithms is considered in
\cite{BirQ88}, but this approach involves significant data movement
and does not scale particularly well.  In \cite{FraGV00}, the
second-stage problems \eqnok{second-stage-lp} are solved concurrently
and inexactly by using an interior-point code. The master process
maintains an upper bound on the optimal objective, and this bound
along with the subgradients obtained from the second-stage problems
yields a polygon whose (approximate) analytic center is calculated
periodically to obtain a new candidate $x$. The approach is based in
part on an algorithm described by Gondzio and Vial~\cite{GonV00}. The
numerical results in \cite{FraGV00} report solution of a two-stage
stochastic linear program with $2.6$ million variables and $1.2$
million constraints in three hours on a cluster of 10 Linux PCs.

\section{L-Shaped Methods} \labtag{sec:lshaped}

We now describe the L-shaped method, a fundamental algorithm for
solving \eqnok{2stage.pl}, and an asynchronous variant.

\subsection{The Multicut L-Shaped Method} \labtag{sec:lshaped:multicut}

The L-shaped method of Van Slyke and Wets~\cite{VanW69} for solving
\eqnok{2stage.pl} proceeds by finding subgradients of partial sums of
the terms that make up $\cQ$ \eqnok{def.Q}, together with linear
inequalities that define the domain of $\cQ$.  The method is
essentially Benders decomposition~\cite{Ben62}, enhanced to deal with
infeasible iterates.  A full description is given in Chapter 5 of
Birge and Louveaux~\cite{BirL97}. We sketch the approach here and
show how it can be implemented in an asynchronous fashion.

We suppose that the second-stage scenarios indexed by $1,2,\dots, N$
are partitioned into $T$ clusters denoted by $\cN_1, \cN_2, \dots,
\cN_T$.  Let $\cQ_{[j]}$ represent the partial sum
from \eqnok{def.Q} corresponding to the cluster $\cN_j$:
\beq \labtag{thetaj}
\cQ_{[j]}(x) = \sum_{i \in \cN_j} p_i \cQ_i(x).
\eeq
The algorithm maintains a model function $m^k_{[j]}$, which is a
piecewise linear lower bound on $\cQ_{[j]}$ for each $j$. We define
this function at iteration $k$ by
\beq \labtag{Qjk}
m_{[j]}^k(x) = \inf \{ \theta_j \, | \, 
 \theta_j e \ge F_{[j]}^k x + f_{[j]}^k \},
\eeq
where $F_{[j]}^k$ is a matrix whose rows are subgradients of
$\cQ_{[j]}$ at previous iterates of the algorithm, and
$e=(1,1,\dots,1)^T$.  The rows of $\theta_j e \ge F_{[j]}^k x +
f_{[j]}^k$ are referred to as {\em optimality cuts}. Upon evaluating
$\cQ_{[j]}$ at the new iterate $x^k$ by solving
\eqnok{second-stage-lp} for each $i \in \cN_j$, a subgradient 
$g_j \in \partial \cQ_{[j]}$ can be obtained from a formula
derived from \eqnok{subg.Qi} and \eqnok{subg.Q}, namely,
\beq \labtag{subg.Qj}
g_j = - \sum_{i \in \cN_j} p_i T(\omega_i)^T \pi(\omega_i),
\eeq
where each $\pi(\omega_i)$ is an optimal dual solution of
\eqnok{second-stage-lp}. 
Since by the subgradient property we have
\[
\cQ_{[j]}(x) \ge g_j^T x + (\cQ_{[j]}(x^k) - g_j^T x^k),
\]
we can obtain $F_{[j]}^{k+1}$ from $F_{[j]}^k$ by appending the row
$g_j^T$, and $f_{[j]}^{k+1}$ from $f_{[j]}^k$ by appending the element
$(\cQ_{[j]}(x^k) - g_j^T x^k)$. In order to keep the number of cuts reasonable,
the cut is not added if $m^k_{[j]}$ is not greater than the value
predicted by the lower bounding approximation (see \eqnok{master}
below).  In this case, the current set of cuts in $F_{[j]}^k$,
$f_{[j]}^k$ adequately models $\cQ_{[j]}$. In addition, we may also
wish to delete some rows from $F_{[j]}^{k+1}$, $f_{[j]}^{k+1}$
corresponding to facets of the epigraph of \eqnok{Qjk} that we do not
expect to be active in later iterations.

The algorithm also maintains a collection of {\em feasibility cuts}
of the form
\beq \labtag{feas.cuts}
D^k x \ge d^k,
\eeq
which have the effect of excluding values of $x$ that were found to be
infeasible, in the sense that some of the second-stage linear programs
\eqnok{second-stage-lp} are infeasible for these values of $x$.  By
Farkas's theorem (see Mangasarian~\cite[p.~31]{Man69}), if the
constraints \eqnok{second-stage-lp.2} are infeasible, there exists
$\pi(\omega_i)$ with the following properties:
\[
W^T \pi(\omega_i) \le 0, \sgap
\left[ h(\omega_i) - T(\omega_i) x \right]^T \pi(\omega_i) > 0.
\]
(In fact, such a $\pi(\omega_i)$ can be obtained from the dual simplex
method for the feasibility problem \eqnok{second-stage-lp.2}.) To
exclude this $x$ from further consideration, we simply add the
inequality $[h(\omega_i) - T(\omega_i) x]^T \pi(\omega_i) \le 0$ to
the constraint set, by appending the row vector $\pi(\omega_i)^T
T(\omega_i)$ to $D^k$ and the element $\pi(\omega_i)^T h(\omega_i)$ to
$d^k$ in \eqnok{feas.cuts}.

The iterate $x^k$ of the multicut L-shaped method is obtained by solving
the following approximation to \eqnok{2stage.pl}:
\beq \labtag{2stage.pl.L}
\min_x \, m_k(x), 
\;\; \mbox{subject to} \; D^k x \ge d^k, \; Ax=b, \; x \ge 0,
\eeq
where 
\beq \labtag{def.mk}
m_k(x) \defeq c^Tx + \sum_{j=1}^T m_{[j]}^k(x).
\eeq
In practice, we substitute from
\eqnok{Qjk} to obtain the following linear program:
 \begin{subequations} \labtag{master}
\beqa
\labtag{master.1}
\min_{x, \theta_1, \dots, \theta_T} \, c^Tx + \sum_{j=1}^T \theta_j, && 
\mbox{subject to} \\
\labtag{master.4}
\theta_j e & \ge & F_{[j]}^k x + f_{[j]}^k, \sgap j=1,2,\dots,T, \\
\labtag{master.3}
D^k x & \ge & d^k, \\
\labtag{master.2}
Ax=b, \;\; x & \ge & 0.
\eeqa
\end{subequations}

The L-shaped method proceeds by solving \eqnok{master} to generate a
new candidate $x$, then evaluating the partial sums \eqnok{thetaj} and
adding optimality and feasibility cuts as described above. The process
is repeated, terminating when the improvement in objective promised by
the subproblem \eqnok{2stage.pl.L} becomes small.

For simplicity we make the following assumption for the remainder of
the paper.
\begin{assumption} \labtag{ass:S}
\mbox{}
\begin{itemize}
\item[(i)] The problem has complete recourse; that is, the feasible
  set of \eqnok{second-stage-lp} is nonempty for all $i=1,2,\dots,N$
  and all $x$, so that the domain of $\cQ(x)$ in \eqnok{def.Q} is $\R^n$.
\item[(ii)] The solution set $\cS$ is nonempty.
\end{itemize}
\end{assumption}
Under this assumption, feasibility cuts of the form \eqnok{feas.cuts},
\eqnok{master.3} do not appear during the course of the algorithm. Our
algorithms and their analysis can be generalized to handle situations
in which Assumption~\ref{ass:S} does not hold, but since our
development is complex enough already, we postpone discussion of these
generalizations to a future report.

Using Assumption~\ref{ass:S}, we can specify the L-shaped algorithm
formally as follows:
\btab
\> {\bf Algorithm LS} \\
\> choose tolerance $\epstol$; \\
\> choose starting point $x^0$; \\
\> define initial model $m_0$ to be a piecewise linear 
underestimate of $\cQ(x)$ \\
\>\>  such that $m_0(x^0) = \cQ(x^0)$ and $m_0$ is bounded below; \\
\> $\cQ_{\rm min} \leftarrow \cQ(x^0)$; \\
\> {\bf for} $k=0,1,2,\dots$ \\
\>\> obtain $x^{k+1}$ by solving \eqnok{2stage.pl.L}; \\
\>\> {\bf if} 
      $\cQ_{\rm min}  - m_k(x^{k+1}) \le \epstol (1+|\cQ_{\rm min}|) $ \\
\>\>\> STOP; \\
\>\> evaluate function and subgradient information at $x^{k+1}$; \\
\>\> $\cQ_{\rm min} \leftarrow  \min(\cQ_{\rm min}, \cQ(x^{k+1}))$; \\
\>\> obtain $m_{k+1}$ by adding optimality cuts to $m_k$; \\
\> {\bf end(for).}
\etab

\subsection{An Asynchronous Parallel Variant of the L-Shaped Method} 
\labtag{sec:lshaped:async}

The L-shaped approach lends itself naturally to implementation in a
master-worker framework. The problem \eqnok{master} is
solved by the master process, while solution of each cluster
$\cN_j$ of second-stage problems, and generation of the associated
cuts, can be carried out by the worker processes running in parallel.
This approach can be adapted for an asynchronous, unreliable
environment in which the results from some second-stage clusters are
not returned in a timely fashion. Rather than having all the worker
processors sit idle while waiting for the tardy results, we can
proceed without them, re-solving the master by using the additional cuts
that were generated by the other second-stage clusters.

We denote the model function simply by $m$ for the asynchronous
algorithm, rather than appending a subscript. Whenever the time comes
to generate a new iterate, the current model is used. In practice, we
would expect the algorithm to give different results each time it is
executed, because of the unpredictable speed and order in which the
functions are evaluated and subgradients generated. Because of
Assumption~\ref{ass:S}, we can write the subproblem 
\beq \labtag{als.subprob}
\min_x \, m(x), 
\;\; \mbox{subject to} \; Ax=b, \; x \ge 0.
\eeq

Algorithm ALS, the asynchronous variant of the L-shaped method that we
describe here, is made up of four key operations, three of which
execute on the master processor and one of which runs on the
workers. These operations are as follows:
\bi
\item {\tt partial\_evaluate}. This is the routine for evaluating
  $\cQ_{[j]}(x)$ defined by \eqnok{thetaj} for a given $x$ and $j$,
  in the process generating a subgradient $g_j$ of $\cQ_{[j]}(x)$. It runs on a
  worker processor and returns its results to the master by
  activating the routine {\tt act\_on\_completed\_task} on the master
  processor.
  
\item {\tt evaluate}. This routine, which runs on the master, simply
  places $T$ tasks of the type {\tt partial\_evaluate} for a given $x$ into the task
  pool for distribution to the worker processors as they become
  available.  The completion of these $T$ tasks is equivalent to evaluating $\cQ(x)$.
  
\item {\tt initialize}. This routine runs on the master processor
  and performs initial bookkeeping, culminating in a call to {\tt
    evaluate} for the initial point $x^0$.

\item {\tt act\_on\_completed\_task}. This routine, which runs on the
master, is activated whenever the results become available from a {\tt
partial\_evaluate} task. It updates the model and increments a counter
to keep track of the number of clusters that have been evaluated at
each candidate point. When appropriate, it solves the master problem
with the latest model to obtain a new candidate iterate\, and will call {\tt evaluate}. 

\ei

In our implementation of both this algorithm and its more
sophisticated cousin Algorithm ATR of Section~\ref{sec:atr}, we may
define a single task to consist of the evaluation of more than one
cluster $\cN_j$. We may bundle, say, $5$ or $10$ clusters into a
single task, in the interests of making the task large enough to
justify the master's effort in packing its data and unpacking its
results, and to maintain the ratio of compute time to communication
cost at a high level. For purposes of simplicity, however, we assume
in the descriptions both of this algorithm and of ATR that each task
consists of a single cluster.

The implementation depends on a ``synchronicity'' parameter $\sigma$
which is the proportion of clusters that must be evaluated at a point
to trigger the generation of a new candidate iterate. Typical values
of $\sigma$ are in the range $0.25$ to $0.9$. A logical variable
${\tt speceval}_k$ keeps track of whether $x^k$ has yet triggered a
new candidate. Initially, ${\tt speceval}_k$ is  set to ${\tt false}$,
then set to ${\tt true}$ when the proportion of evaluated clusters
passes the threshold $\sigma$.  

We now specify all the methods making up Algorithm ALS.

\btab
\>{\bf ALS:} \ \ {\tt  partial\_evaluate}$(x^q,q,j,\cQ_{[j]}(x^q),g_j)$ \\
\> Given $x^q$, index  $q$, and  partition number $j$, 
evaluate $\cQ{[j]}(x^q)$ from \eqnok{thetaj} \\
\>\> together with a partial subgradient $g_j$ from \eqnok{subg.Qj};
\\
\> Activate {\tt act\_on\_completed\_task}$(x^q,q,j,\cQ_{[j]}(x^q),g_j)$
on the master processor.
\etab

\medskip

\btab
\> {\bf ALS:} \ \ {\tt  evaluate}$(x^q,q)$ \\
\> {\bf for} $j=1,2,\dots, T$ (possibly concurrently) \\
\>\> {\tt partial\_evaluate}$(x^q,q,j,\cQ_{[j]}(x^q), g_j)$; \\
\> {\bf end (for)}
\etab

\medskip

\btab
\> {\bf ALS:} \ \ {\tt initialize} \\
\> choose tolerance $\epstol$; \\
\> choose starting point $x^0$; \\
\> choose threshold $\sigma \in (0,1]$; \\
\> $\cQ_{\rm min} \leftarrow \infty$; \\
\> $k \leftarrow 0$, ${\tt speceval}_0 \leftarrow {\tt false}$, $t_0 \leftarrow 0$; \\
\> {\tt evaluate}$(x^0,0)$.
\etab

\medskip

\btab
\> {\bf ALS:} \ \  
{\tt act\_on\_completed\_task}$(x^q,q,j,\cQ_{[j]}(x^q),g_j)$ \\
\> $t_q \leftarrow t_q+1$; \\
\> add $\cQ_{[j]}(x^q)$ and cut $g_j$ to the model $m$; \\
\> {\bf if}  $t_q = T$ \\
\>\> $\cQ_{\rm min} \leftarrow \min ( \cQ_{\rm min}, \cQ(x^q))$; \\
\> {\bf else if}  $t_q \ge \sigma T$ {\bf and} not ${\tt speceval}_q$ \\
\>\> ${\tt speceval}_q \leftarrow ${\tt true}; \\
\>\> $k \leftarrow k+1$;  \\
\>\> solve  current model problem \eqnok{als.subprob} to obtain $x^{k+1}$; \\
\>\> {\bf if} $\cQ_{\rm min}  - m(x^{k+1}) \le \epstol (1+|\cQ_{\rm min}|) $ \\
\>\>\> STOP; \\
\>\> {\tt evaluate}$(x^k,k)$; \\
\> {\bf end (if)}

\etab

We present results for Algorithm ALS in Section~\ref{sec:results}.
While the algorithm is able to use a large number of worker processors
on our opportunistic platform, it suffers from the usual drawbacks of
the L-shaped method, namely, that cuts, once generated, must be
retained for the remainder of the computation to ensure convergence
and that large steps are typically taken on early iterations before a
sufficiently good model approximation to $\cQ(x)$ is created, making
it impossible to exploit prior knowledge about the location of the
solution.

\section{A Bundle-Trust-Region Method} \labtag{sec:tr}

Trust-region approaches can be implemented by making only minor
modifications to implementations of the L-shaped method, and they
possesses several practical advantages along with stronger convergence
properties. The trust-region methods we describe here are related to
the regularized decomposition method of Ruszczy{\'n}ski~\cite{Rus86}
and the bundle-trust-region approaches of Kiwiel~\cite{Kiw90} and
Hirart-Urruty and Lemar\'echal~\cite[Chapter~XV]{HirL93}. The main
differences are that we use box-shaped trust regions yielding linear
programming subproblems (rather than quadratic programs) and that our
methods manipulate the size of the trust region directly rather than
indirectly via a regularization parameter.

When requesting a subgradient of $\cQ$ at some
point $x$, our algorithms do not require particular (e.g., extreme)
elements of the subdifferential to be supplied. Nor do they require
the subdifferential $\partial \cQ(x)$ to be representable as a convex
combination of a finite number of vectors. In this respect, our
algorithms contrast with that of Ruszczy{\'n}ski~\cite{Rus86}, for
instance, which exploits the piecewise-linear nature of the objectives
$\cQ_i$ in \eqnok{second-stage-lp}. Because of our weaker conditions
on the subgradient information, we cannot prove a finite termination
result of the type presented in \cite[Section~3]{Rus86}.  However,
these conditions potentially allow our algorithms to be extended to a
more general class of convex nondifferentiable functions. We hope to
explore these generalizations in future work.

\subsection{A Method Based on $\ell_{\infty}$ Trust Regions}
\labtag{sec:tr:tr}

A key difference between the trust-region approach of this section and
the L-shaped method of the preceding section is that we impose an
$\ell_{\infty}$ norm bound on the size of the step. It is implemented
by simply adding bound constraints to the linear programming
subproblem \eqnok{master} as follows:
\beq \labtag{master.tr.bounds}
-\Delta e \le x-x^k \le \Delta e,
\eeq
where $e=(1,1,\dots,1)^T$, $\Delta$ is the trust-region radius, and
$x^k$ is the current iterate. During the $k$th iteration, it may be
necessary to solve several problems with trust regions of the form
\eqnok{master.tr.bounds}, with different model functions $m$ and
possibly different values of $\Delta$, before a satisfactory new
iterate $x^{k+1}$ is identified. We refer to $x^k$ and $x^{k+1}$ as
{\em major iterates} and the points $x^{k,\ell}$, $\ell=0,1,2,\dots$
obtained by minimizing the current model function subject to the
constraints and trust-region bounds of the form
\eqnok{master.tr.bounds} as {\em minor iterates}. Another key
difference between the trust-region approach and the L-shaped approach
is that a minor iterate $x^{k,\ell}$ is accepted as the new major
iterate $x^{k+1}$ only if it yields a substantial reduction in the
objective function $\cQ$ over the previous iterate $x^k$, in a sense
to be defined below.  A further important difference is that one can
delete optimality cuts from the model functions, between minor and
major iterations, without compromising the convergence properties of
the algorithm.

To specify the method, we need to augment the notation established in
the previous section.  We define $m_{k,\ell}(x)$ to be the model
function after $\ell$ minor iterations have been performed at
iteration $k$, and $\Delta_{k,\ell}>0$ to be the trust-region radius
at the same stage.  Under Assumption~\ref{ass:S}, there are no
feasibility cuts, so that the problem to be solved to obtain the minor
iteration $x^{k,\ell}$ is as follows:
\beq \labtag{trsub.kl}
\min_x \, m_{k,\ell}(x) \;\; \mbox{subject to} \;Ax=b, \; x \ge 0, \; 
\| x-x^k \|_{\infty} \le \Delta_{k,\ell}
\eeq
(cf. \eqnok{2stage.pl.L}). By expanding this problem in a similar
fashion to \eqnok{master}, we obtain
\begin{subequations} \labtag{master.kl}
\beqa
\labtag{master.kl.1}
\min_{x, \theta_1, \dots, \theta_T} \, c^Tx + \sum_{j=1}^T \theta_j, && 
\mbox{subject to} \\
\labtag{master.kl.4}
\theta_j e & \ge & F_{[j]}^{k,\ell} x + f_{[j]}^{k,\ell}, \sgap j=1,2,\dots,T, \\
\labtag{master.kl.2}
Ax=b, \;\; x & \ge & 0, \\
\labtag{master.kl.tr}
-\Delta_{k,\ell} e \le  x-x^k & \le & \Delta_{k,\ell} e.
\eeqa
\end{subequations}

We assume the initial model $m_{k,0}$ at major iteration $k$ to
satisfy the following two properties:
\begin{subequations} \labtag{mkprop}
\beqa \labtag{mkprop.1}
& m_{k,0}(x^k) = \cQ(x^k), \\
\labtag{mkprop.2}
& \mbox{$m_{k,0}$ is a piecewise linear underestimate of $\cQ$}.
\eeqa
\end{subequations}

Denoting the solution of the subproblem \eqnok{master.kl} by
$x^{k,\ell}$, we accept this point as the new iterate $x^{k+1}$ if the
decrease in the actual objective $\cQ$ (see \eqnok{2stage.pl}) is at
least some fraction of the decrease predicted by the model function
$m_{k,\ell}$. That is, for some constant $\xi \in (0,1/2)$, the
acceptance test is
\beq \labtag{tr.accept}
\cQ(x^{k,\ell}) \le \cQ(x^k) - \xi 
\left( \cQ(x^k) - m_{k,\ell}(x^{k,\ell}) \right).
\eeq
(A typical value for $\xi$ is $10^{-4}$.)  

If the test \eqnok{tr.accept} fails to hold, we obtain a new model
function $m_{k,\ell+1}$ by adding and possibly deleting cuts from
$m_{k,\ell}(x)$. This process aims to refine the model function, so
that it eventually generates a new major iteration, while economizing
on storage by allowing deletion of subgradients that no longer seem
helpful. Addition and deletion of cuts are implemented by adding and
deleting rows from $F_{[j]}^{k,\ell}$ and $f_{[j]}^{k,\ell}$, to
obtain $F_{[j]}^{k,\ell+1}$ and $f_{[j]}^{k,\ell+1}$, for
$j=1,2,\dots,T$.

Given some parameter $\eta \in [0,1)$, we obtain $m_{k,\ell+1}$ from
$m_{k,\ell}$ by means of the following procedure:
\btab
\> {\bf Procedure Model-Update} $(k,\ell)$ \\
\> {\bf for each} optimality cut \\
\>\> {\tt possible\_delete}  $\leftarrow$ {\tt true}; \\
\>\> {\bf if} the cut was generated at $x^k$ \\
\>\>\> {\tt possible\_delete}  $\leftarrow$ {\tt false}; \\
\>\> {\bf else if} the cut is active at the solution of \eqnok{master.kl} \\
\>\>\> {\tt possible\_delete}  $\leftarrow$ {\tt false}; \\
\>\> {\bf else if} the cut was generated at an earlier minor iteration \\
\>\>\>
$\bar{\ell}=0,1,\dots,\ell-1$ such that
\etab
\beq \labtag{cut.delete.criterion}
\cQ(x^k) - m_{k,\ell}(x^{k,\ell}) > \eta 
\left[ \cQ(x^k) - m_{k,\bar\ell}(x^{k,\bar\ell}) \right]
\eeq
\btab
\>\>\> {\tt possible\_delete}  $\leftarrow$ {\tt false}; \\
\>\> {\bf end (if)} \\
\>\> {\bf if} {\tt possible\_delete} \\
\>\>\> possibly delete the cut; \\
\> {\bf end (for each)} \\
\> add optimality cuts obtained from each of the component functions \\
\>\> $\cQ_{[j]}$ at $x^{k,\ell}$. \\
\etab

In our implementation, we delete the cut if ${\tt possible\_delete}$
is true at the final conditional statement and, in addition, the cut
has not been active during the last 100 solutions of
\eqnok{master.kl}. More details are given in
Section~\ref{sec:results:parameters}.

Because we retain all cuts active at $x^k$ during the course of
major iteration $k$, the following extension of \eqnok{mkprop.1} holds:
\beq \labtag{mkprop.1a}
m_{k,\ell}(x^k) = \cQ(x^k), \;\; \ell=0,1,2,\dots.
\eeq
Since we add only subgradient information, the following
generalization of \eqnok{mkprop.2} also holds uniformly:
\beq \labtag{mkprop.2a}
\mbox{$m_{k,\ell}$ is a piecewise linear underestimate of $\cQ$, for $\ell=0,1,2,\dots.$}
\eeq

We may also decrease the trust-region radius $\Delta_{k,\ell}$ between
minor iterations (that is, choose $\Delta_{k,\ell+1} <
\Delta_{k,\ell}$) when the test \eqnok{tr.accept} fails to hold. We do
so if the match between model and objective appears to be particularly
poor.  If $\cQ(x^{k,\ell})$ exceeds $\cQ(x^k)$ by more than an
estimate of the quantity
\beq \labtag{reduce.delta.1}
\max_{\| x-x^k\|_{\infty} \le 1} \, \cQ(x^k) - \cQ(x),
\eeq
we conclude that the ``upside'' variation of the function $\cQ$
deviates too much from its ``downside'' variation, and we choose the
new radius $\Delta_{k,\ell+1}$ to bring these quantities more nearly
into line. Our estimate of \eqnok{reduce.delta.1} is simply
\[
\frac{1}{\min(1,\Delta_{k,\ell})} 
\left[ \cQ(x^k)  - m_{k,\ell}(x^{k,\ell}) \right],
\]
that is, an extrapolation of the model reduction on the current trust
region to a trust region of radius $1$.  Our complete strategy for
reducing $\Delta$ is therefore as follows. (The {\tt counter} is
initialized to zero at the start of each major iteration.)
\btab
\> {\bf Procedure Reduce-$\Delta$} \\
\> evaluate
\etab
\beq \labtag{reduce.delta.2}
\rho = {\min(1,\Delta_{k,\ell})} \frac{\cQ(x^{k,\ell}) - \cQ(x^k)}{\cQ(x^k)  - m_{k,\ell}(x^{k,\ell})};
\eeq
\btab
\> {\bf if} $\rho>0$ \\
\>\> {\tt counter} $\leftarrow$ {\tt counter}$+1$; \\
\> {\bf if} $\rho>3$ {\bf or} 
({\tt counter} $\ge 3$ {\bf and} $\rho \in (1,3]$) \\
\>\> set
\etab
\[
\Delta_{k,\ell+1} = \frac{1}{\min(\rho,4)} \Delta_{k,\ell};
\]
\btab
\>\> reset {\tt counter} $\leftarrow 0$; 
\etab
This procedure is related to the technique of
Kiwiel~\cite[p.~109]{Kiw90} for increasing the coefficient of the
quadratic penalty term in his regularized bundle method.

If the test \eqnok{tr.accept} is passed, so that we have
$x^{k+1} = x^{k,\ell}$, we have a great deal of flexibility in
defining the new model function $m_{k+1,0}$.  We require only that the
properties \eqnok{mkprop} are satisfied, with $k+1$ replacing $k$.
Hence, we are free to delete much of the optimality cut information
accumulated at iteration $k$ (and previous iterates). In practice, of
course, it is wise to delete only those cuts that have been inactive
for a substantial number of iterations; otherwise we run the risk that
many new function and subgradient evaluations will be required to
restore useful model information that was deleted prematurely.

If the step to the new major iteration $x^{k+1}$ shows a particularly
close match between the true function $\cQ$ and the model function
$m_{k,\ell}$ at the last minor iteration of iteration $k$, we consider
increasing the trust-region radius. Specifically, if
\beq \labtag{tr.incr.1}
\cQ(x^{k,\ell}) \le \cQ(x^k) - 0.5
\left( \cQ(x^k) - m_{k,\ell}(x^{k,\ell}) \right), \sgap
\| x^k - x^{k,\ell} \|_{\infty} = \Delta_{k,\ell},
\eeq
then we set
\beq \labtag{tr.incr.3}
\Delta_{k+1,0} = \min ( \Delta_{\rm hi}, 2 \Delta_{k,\ell}),
\eeq
where $\Delta_{\rm hi}$ is a prespecified upper bound on the radius.

Before specifying the algorithm formally, we define the convergence
test. Given a parameter $\epstol>0$, we terminate if 
\beq \labtag{conv.test}
\cQ(x^k) - m_{k,\ell}(x^{k,\ell}) \le 
\epstol  (1+ |\cQ(x^k)|).
\eeq

\btab
\> {\bf Algorithm TR} \\
\> choose $\xi \in (0,1/2)$, maximum trust region $\Delta_{\rm hi}$, 
tolerance $\epstol$; \\
\> choose starting point $x^0$; \\
\> define initial model $m_{0,0}$ with the properties \eqnok{mkprop} (for $k=0$); \\
\> choose $\Delta_{0,0} \in (0, \Delta_{\rm hi}]$; \\
\> {\bf for} $k=0,1,2,\dots$ \\
\>\> {\tt finishedMinorIteration} $\leftarrow$ {\tt false}; \\
\>\> $\ell \leftarrow 0$; ${\tt counter} \leftarrow 0$; \\
\>\> {\bf repeat} \\
\>\>\> solve \eqnok{trsub.kl} to obtain $x^{k,\ell}$; \\
\>\>\> {\bf if} \eqnok{conv.test} is satisfied \\
\>\>\>\> STOP with approximate solution $x^k$; \\
\>\>\> evaluate function and subgradient at $x^{k,\ell}$; \\
\>\>\> {\bf if} \eqnok{tr.accept} is satisfied \\
\>\>\>\> set  $x^{k+1} = x^{k,\ell}$; \\
\>\>\>\> obtain $m_{k+1,0}$ by possibly deleting cuts from $m_{k,\ell}$, but \\
\>\>\>\>\> 
 retaining the properties \eqnok{mkprop} (with $k+1$ replacing $k$); \\
\>\>\>\> choose $\Delta_{k+1,0} \in [ \Delta_{k,\ell}, \Delta_{\rm hi}]$ 
according to \eqnok{tr.incr.1}, \eqnok{tr.incr.3}; \\
\>\>\>\> {\tt finishedMinorIteration} $\leftarrow$ {\tt true}; \\
\>\>\> {\bf else} \\
\>\>\>\> obtain  $m_{k,\ell+1}$ from  $m_{k,\ell}$ 
via Procedure Model-Update $(k,\ell)$; \\
\>\>\>\> obtain $\Delta_{k,\ell+1}$ via Procedure Reduce-$\Delta$; \\
\>\>\> $\ell \leftarrow \ell+1$; \\
\>\> {\bf until} {\tt finishedMinorIteration} \\
\> {\bf end (for)}
\etab

\subsection{Analysis of the Trust-Region Method}
\labtag{sec:tr:analysis}

We now describe the convergence properties of Algorithm TR. We show
that for $\epstol=0$, the algorithm either terminates at a solution
or generates a sequence of major iterates that approaches the
solution set $\cS$ (Theorem~\ref{th:tr:conv}). When $\epstol > 0$, the
algorithm terminates finitely; that is, it avoids generating infinite
sequences either of major or minor iterates (Theorem~\ref{th:fint}).

Given some starting point $x^0$ satisfying the constraints
$Ax^0=b$, $x^0 \ge 0$, and setting $\cQ_0 = \cQ(x^0)$, we define the
following quantities that are useful in describing and analyzing the
algorithm:
\beqa
\labtag{def.ls}
\cL(\cQ_0) &=& \{ x \, | \, Ax=b, x \ge 0, \cQ(x) \le \cQ_0 \}, \\
\labtag{def.lsn}
\cL(\cQ_0;\Delta) &=& \{ x \, | \, \|x-y \| \le \Delta, \, 
\mbox{for some  $y \in \cL(\cQ_0)$} \}, \\
\labtag{def.beta} 
\beta &=& \sup \{ \| g \|_1 \, | \, g \in \partial \cQ(x), \, 
\mbox{for some $x \in \cL(\cQ_0;\Delta_{\rm hi})$} \}.
\eeqa
Using Assumption~\ref{ass:S}, we can easily show that $\beta <
\infty$.

We start by showing that the optimal objective value for
\eqnok{trsub.kl} cannot decrease from one minor iteration to the next.
\begin{lemma} \labtag{lem:mkl}
  Suppose that $x^{k,\ell}$ does not satisfy the acceptance test
  \eqnok{tr.accept}. Then we have
\[
m_{k,\ell}(x^{k,\ell}) \le m_{k,\ell+1}(x^{k,\ell+1}).
\]
\end{lemma}
\begin{proof}
  In obtaining $m_{k,\ell+1}$ from $m_{k,\ell}$ in Model-Update, we do
  not allow deletion of cuts that were active at the solution
  $x^{k,\ell}$ of \eqnok{master.kl}. Using $\bar{F}_{[j]}^{k,\ell}$
  and $\bar{f}_{[j]}^{k,\ell}$ to denote the active rows in
  $F_{[j]}^{k,\ell}$ and $f_{[j]}^{k,\ell}$, we have that $x^{k,\ell}$
  is also the solution of the following linear program (in which the
  inactive cuts are not present):
\begin{subequations} \labtag{master2.kl}
\beqa
\labtag{master2.kl.1}
\min_{x, \theta_1, \dots, \theta_T} \, c^Tx + \sum_{j=1}^T \theta_j, && 
\mbox{subject to} \\
\labtag{master2.kl.4}
\theta_j e & \ge & \bar{F}_{[j]}^{k,\ell} x + \bar{f}_{[j]}^{k,\ell}, \sgap j=1,2,\dots,T, \\
\labtag{master2.kl.2}
Ax=b, \;\; x & \ge & 0, \\
\labtag{master2.kl.tr}
-\Delta_{k,\ell} e \le  x-x^k & \le & \Delta_{k,\ell} e.
\eeqa
\end{subequations}
The subproblem to be solved for $x^{k,\ell+1}$ differs from
\eqnok{master2.kl} in two ways. First, additional rows may be added to
$\bar{F}_{[j]}^{k,\ell}$ and $\bar{f}_{[j]}^{k,\ell}$, consisting of
function values and subgradients obtained at $x^{k,\ell}$ and also
inactive cuts carried over from the previous \eqnok{master.kl}. Second,
the trust-region radius $\Delta_{k,\ell+1}$ may be smaller than
$\Delta_{k,\ell}$. Hence, the feasible region of the problem to be
solved for $x^{k,\ell+1}$ is a subset of the feasible region for
\eqnok{master2.kl}, so the optimal objective value cannot be smaller. 
\end{proof}

Next we have a result about the amount of reduction in the model
function $m_{k,\ell}$.
\begin{lemma} \labtag{lem:tr:1}
For all $k=0,1,2,\ldots$ and $\ell=0,1,2,\ldots$, we have that 
\begin{subequations} \labtag{lem:tr:inequalities}
\beqa 
\nonumber
m_{k,\ell}(x^k) - m_{k,\ell}(x^{k,\ell}) &=&
\cQ(x^k) - m_{k,\ell}(x^{k,\ell}) \\
\labtag{tr.2a}
& \ge  &
\min \left( \Delta_{k,\ell}, \| x^k - P(x^k)\|_{\infty} \right)
\frac{\cQ(x^k) - \cQ^*}{\| x^k - P(x^k) \|_{\infty}} \\
\labtag{tr.2b}
& \ge &
\hat{\epsilon} \min \left( \Delta_{k,\ell}, \| x^k - P(x^k)\|_{\infty} \right),
\eeqa
\end{subequations}
where $\hat{\epsilon}>0$ is defined in \eqnok{weak.sharp}.
\end{lemma}
\begin{proof}
  The first equality follows immediately from \eqnok{mkprop.1a}, while
  the second inequality \eqnok{tr.2b} follows immediately from
  \eqnok{tr.2a} and \eqnok{weak.sharp}. We now prove \eqnok{tr.2a}.
 
Consider the following subproblem in the scalar $\tau$:
\beq \labtag{tr.3}
\min_{\tau \in [0,1]} \, m_{k,\ell} \left( x^k + \tau [ P(x^k) - x^k] \right)
\;\; \mbox{subject to} \; \left\| \tau [ P(x^k) - x^k] \right\|_{\infty} \le
\Delta_{k,\ell}.
\eeq
Denoting the solution of this problem by $\tau_{k,\ell}$, we have by
comparison with \eqnok{trsub.kl} that
\beq \labtag{tr.3a}
m_{k,\ell} (x^{k,\ell}) \le 
m_{k,\ell} \left( x^k  + \tau_{k,\ell} [ P(x^k) - x^k] \right).
\eeq
If $\tau=1$ is feasible in \eqnok{tr.3}, we have  from \eqnok{tr.3a} and
\eqnok{mkprop.2a} that
\beqas
\lefteqn{m_{k,\ell} (x^{k,\ell}) \le 
m_{k,\ell} \left( x^k  + \tau_{k,\ell} [ P(x^k) - x^k] \right) } \\
& \le &
m_{k,\ell} \left( x^k  + [ P(x^k) - x^k] \right)
= m_{k,\ell} (P(x^k)) \le \cQ(P(x^k)) = \cQ^*.
\eeqas
Therefore, when $\tau=1$ is feasible for \eqnok{tr.3}, we have from
\eqnok{mkprop.1a} that
\[
m_{k,\ell}(x^k) - m_{k, \ell}(x^{k,\ell}) \ge \cQ(x^k) - \cQ^*,
\]
so that  \eqnok{tr.2a} holds in this case.

When $\tau=1$ is infeasible for \eqnok{tr.3}, consider setting $\tau =
\Delta_{k,\ell} / \| x^k-P(x^k) \|_{\infty}$ (which is certainly feasible for
\eqnok{tr.3}). We have from \eqnok{tr.3a}, the definition of
$\tau_{k,\ell}$, the fact \eqnok{mkprop.2a} that $m_{k,\ell}$
underestimates $\cQ$, and convexity of $\cQ$ that
\beqas
m_{k,\ell}(x^{k,\ell})
& \le & m_{k,\ell} 
\left( x^k + \Delta_{k,\ell} \frac{P(x^k)-x^k}{\|P(x^k)-x^k\|_{\infty}} \right) \\
& \le & \cQ
\left( x^k + \Delta_{k,\ell} \frac{P(x^k)-x^k}{\|P(x^k)-x^k\|_{\infty}} \right) \\
& \le & \cQ(x^k) + 
\frac{\Delta_{k,\ell}}{\|P(x^k)-x^k\|_{\infty}} (\cQ^* - \cQ(x^k)).
\eeqas
Therefore, using \eqnok{mkprop.1a}, we have
\[
m_{k,\ell}(x^k) - m_{k,\ell}(x^{k,\ell}) \ge
\frac{\Delta_{k,\ell}}{\|P(x^k)-x^k\|_{\infty}} [ \cQ(x^k) - \cQ^* ],
\]
verifying \eqnok{tr.2a} in this case as well.
\end{proof}

Our next result finds a lower bound on the trust-region radii
$\Delta_{k,\ell}$. For purposes of this result we define a quantity
$E_k$ to measure the closest approach to the solution set for all
iterates up to and including $x^k$, that is,
\beq \labtag{def:Ek}
E_k \defeq \min_{\bar{k}=0,1,\dots,k} 
\| x^{\bar{k}} - P(x^{\bar{k}}) \|_{\infty}.
\eeq
Note that $E_k$ decreases monotonically with $k$. We also define
$\Delta_{\rm init}$ to be the initial value of the trust region.
\begin{lemma} \labtag{lem:trbounds}
  There is a constant $\Delta_{\rm lo} >0$ such that for all trust
  regions $\Delta_{k,\ell}$ used in the course of Algorithm TR, we
  have
\[
\Delta_{k,\ell} \ge \min( \Delta_{\rm lo}, E_k/4).
\]
\end{lemma}
\begin{proof}
  We prove the result by showing that the value $\Delta_{\rm lo} =
  (1/4) \min(1, \Delta_{\rm init}, \hat{\epsilon}/\beta)$ has the
  desired property, where $\hat{\epsilon}$ is from \eqnok{weak.sharp}
  and $\beta$ is from \eqnok{def.beta}.
  
  Suppose for contradiction that there are indices $k$ and $\ell$ such
  that
\[
\Delta_{k, \ell} < \frac14 \min  
\left( 1, \frac{\hat{\epsilon}}{\beta}, \Delta_{\rm init}, E_k \right).
\]
Since the trust region can be reduced by at most a factor of $4$ by
Procedure Reduce-$\Delta$, there must be an earlier trust region
radius $\Delta_{\bar{k}, \bar{\ell}}$ (with $\bar{k} \le k$) such that
\beq \labtag{poo.3}
\Delta_{\bar{k},\bar{\ell}} < 
\min \left( 1, \frac{\hat{\epsilon}}{\beta}, E_{k} \right),
\eeq
and $\rho>1$ in \eqnok{reduce.delta.2}, that is,
\beqa
\nonumber 
\cQ(x^{\bar{k},\bar{\ell}}) - \cQ(x^{\bar{k}}) & > & 
\frac{1}{\min(1,\Delta_{\bar{k},\bar{\ell}})}
\left( \cQ(x^{\bar{k}}) - 
m_{\bar{k},\bar{\ell}}(x^{\bar{k},\bar{\ell}}) \right) \\
\labtag{poo.4}
& = & \frac{1}{\Delta_{\bar{k},\bar{\ell}}}
\left( 
\cQ(x^{\bar{k}}) - m_{\bar{k},\bar{\ell}}(x^{\bar{k},\bar{\ell}}) 
\right).
\eeqa
By applying Lemma~\ref{lem:tr:1}, and using  \eqnok{poo.3}, we have
\beq \labtag{poo.4a}
\cQ(x^{\bar{k}}) - m_{\bar{k},\bar{\ell}}(x^{\bar{k},\bar{\ell}}) \ge
\hat{\epsilon} \min \left( \Delta_{\bar{k},\bar{\ell}}, 
\| x^{\bar{k}} - P(x^{\bar{k}}) \|_{\infty} \right) =
\hat{\epsilon} \Delta_{\bar{k},\bar{\ell}}
\eeq
where the last equality follows from 
$\| x^{\bar{k}} - P(x^{\bar{k}}) \|_{\infty} \ge E_{\bar{k}} \ge E_k$ and
\eqnok{poo.3}.
By combining  \eqnok{poo.4a} with  \eqnok{poo.4}, we have that 
\beq \labtag{poo.5}
\cQ(x^{\bar{k},\bar{\ell}}) - \cQ(x^{\bar{k}}) > \hat{\epsilon}.
\eeq
By using standard properties of subgradients, we have
\beqa 
\nonumber
\lefteqn{\cQ(x^{\bar{k},\bar{\ell}}) - \cQ(x^{\bar{k}}) \le 
g_{\bar{\ell}}^T(x^{\bar{k},\bar{\ell}} - x^{\bar{k}})} \\
\labtag{subd.5}
& \le &
\| g_{\bar{\ell}} \|_1 \| x^{\bar{k}} - x^{\bar{k},\bar{\ell}} \|_{\infty} 
\le \| g_{\bar{\ell}} \|_1 \Delta_{\bar{k},\bar{\ell}}, \;\;
\mbox{for all} \; g_{\bar{\ell}} \in \partial \cQ(x^{\bar{k},\bar{\ell}}).
\eeqa
By combining this expression with \eqnok{poo.5}, and using
\eqnok{poo.3} again, we obtain that
\[
\| g_{\bar{\ell}} \|_1 \ge 
\frac{\hat{\epsilon}}{\Delta_{\bar{k},\bar{\ell}}} > \beta.
\]
However, since $x^{\bar{k},\bar{\ell}} \in \cL(\cQ_0;\Delta_{\rm hi})$, we have
from \eqnok{def.beta} that $\| g_{\bar{\ell}} \|_1 \le \beta$, giving a
contradiction.
\end{proof}

Finite termination of the inner iterations is proved in the following
two results. Recall that the parameters $\xi$ and $\eta$ are defined
in \eqnok{tr.accept} and \eqnok{cut.delete.criterion}, respectively.
\begin{lemma} \labtag{lem:tr:ft}
  Let $\epstol=0$ in Algorithm TR, and let $\bar{\eta}$ be
  any constant satisfying $0<\bar{\eta}<1$, $\bar{\eta}>\xi$,
  $\bar{\eta} \ge \eta$. Let $\ell_1$ be any index such that
  $x^{k,\ell_1}$ fails to satisfy the test \eqnok{tr.accept}.  Then
  either the sequence of inner iterations eventually yields a point
  $x^{k,\ell_2}$ satisfying the acceptance test \eqnok{tr.accept}, or
  there is an index $\ell_2>\ell_1$ such that
\beq \labtag{tr.6}
\cQ(x^k) - m_{k,\ell_2}(x^{k,\ell_2}) \le \bar{\eta} \left[
\cQ(x^k) - m_{k,\ell_1}(x^{k,\ell_1}) \right].
\eeq
\end{lemma}
\begin{proof}
  Suppose for contradiction that the none of the minor iterations
  following $\ell_1$ satisfies either \eqnok{tr.accept} or the
  criterion \eqnok{tr.6}; that is,
\beqa \nonumber 
\cQ(x^k) - m_{k,q}(x^{k,q}) & > & \bar{\eta} \left[
\cQ(x^k) - m_{k,\ell_1}(x^{k,\ell_1}) \right],  \\
\labtag{contra}
& \ge & \eta \left[ \cQ(x^k) - m_{k,\ell_1}(x^{k,\ell_1}) \right],
\;\; \mbox{\rm for all $q > \ell_1$}.
\eeqa
It follows from this bound, together with Lemma~\ref{lem:mkl} and
Procedure Model-Update, that none of the cuts generated at minor
iterations $q \ge \ell_1$ is deleted.  

We assume in the remainder of the proof that $q$ and $\ell$ are
generic minor iteration indices that satisfy
\[
q > \ell \ge \ell_1.
\]

Because the function and subgradients from minor iterations
$x^{k,\ell}$, $l=l_1,l_1+1, \dots$ are retained throughout the major
iteration $k$, we have
\beq \labtag{matchQ}
m_{k,q}(x^{k,\ell}) = \cQ(x^{k,\ell}).
\eeq
By definition of the subgradient, we have 
\beq \labtag{subgrad.mkq}
m_{k,q}(x) - m_{k,q}(x^{k,\ell}) \ge g^T (x-x^{k,\ell}), \;\;
\mbox{for all} \; g \in \partial m_{k,q}(x^{k,\ell}).
\eeq
Therefore, from \eqnok{mkprop.2a} and \eqnok{matchQ}, it follows that 
\[
\cQ(x)-\cQ(x^{k,\ell}) \ge g^T (x-x^{k,\ell}), \;\; \mbox{for all} \;
g \in \partial m_{k,q}(x^{k,\ell}),
\]
so that
\beq \labtag{mkqQ}
\partial m_{k,q}(x^{k,\ell}) \subset \partial \cQ(x^{k,\ell}).
\eeq

Since $\cQ(x^k) < \cQ(x^0) = \cQ_0$, we have from \eqnok{def.ls} that
$x^k \in \cL(\cQ_0)$. Therefore, from the definition \eqnok{def.lsn}
and the fact that $\| x^{k,\ell} - x^k \| \le \Delta_{k,\ell} \le
\Delta_{\rm hi}$, we have that $x^{k,\ell} \in \cL(\cQ_0;\Delta_{\rm
hi})$. It follows from \eqnok{def.beta} and \eqnok{mkqQ} that
\beq \labtag{gbound}
\| g \|_1 \le \beta, \;\; \mbox{for all} \; g \in \partial m_{k,q}(x^{k,\ell}).
\eeq

Since $x^{k,\ell}$ is rejected by the test \eqnok{tr.accept}, we
have from \eqnok{matchQ} and Lemma~\ref{lem:mkl} that the following
inequalities hold:
\beqas
m_{k,q}(x^{k,\ell}) = \cQ(x^{k,\ell}) 
& \ge &\cQ(x^k) - \xi \left[ \cQ(x^k) - m_{k,\ell}(x^{k,\ell}) \right] \\
& \ge & \cQ(x^k) - \xi \left[ \cQ(x^k) - m_{k,\ell_1}(x^{k,\ell_1}) \right].
\eeqas
By rearranging this expression, we obtain
\beq \labtag{tr.8}
\cQ(x^k) - m_{k,q}(x^{k,\ell}) \le 
\xi \left[ \cQ(x^k) - m_{k,\ell_1}(x^{k,\ell_1}) \right].
\eeq

Consider now all points $x$ satisfying 
\beq \labtag{xkl.nbd}
\| x-x^{k,\ell} \|_{\infty} \le
\frac{\bar{\eta}-\xi}{\beta} 
\left[ \cQ(x^k)-m_{k,\ell_1}(x^{k,\ell_1}) \right]
\defeq \zeta>0.
\eeq
Using this bound together with \eqnok{subgrad.mkq} and \eqnok{gbound},
we obtain
\beqas
\lefteqn{ m_{k,q}(x^{k,\ell}) - m_{k,q}(x) \le g^T(x^{k,\ell} - x ) } \\
& \le & \beta \| x^{k,\ell}-x \|_{\infty} 
\le (\bar{\eta} - \xi) \left[ \cQ(x^k)-m_{k,\ell_1}(x^{k,\ell_1}) \right].
\eeqas
By combining this bound with  \eqnok{tr.8}, we find that the following
bound is satisfied for all $x$ in the neighborhood \eqnok{xkl.nbd}:
\beqas
\cQ(x^k) - m_{k,q}(x) &=& 
\left[ \cQ(x^k) - m_{k,q}(x^{k,\ell}) \right] +  
\left[ m_{k,q}(x^{k,\ell}) - m_{k,q}(x) \right] \\
 & \le & \bar{\eta} \left[ \cQ(x^k)-m_{k,\ell_1}(x^{k,\ell_1}) \right].
\eeqas
It follows from this bound, in conjunction with \eqnok{contra}, that
$x^{k,q}$ (the solution of the trust-region problem with model
function $m_{k,q}$) cannot lie in the neighborhood \eqnok{xkl.nbd}.
Therefore, we have
\beq \labtag{meshprop}
\| x^{k,q} - x^{k,\ell} \|_{\infty} > \zeta.
\eeq
But since $\| x^{k,\ell} - x^k \|_{\infty} \le \Delta_k \le
\Delta_{\rm hi}$ for all $\ell \ge \ell_1$, it is impossible for an
infinite sequence $\{ x^{k,\ell} \}_{\ell \ge \ell_1}$ to satisfy
\eqnok{meshprop}. We conclude that \eqnok{tr.6} must hold for some
$\ell_2 \ge \ell_1$, as claimed.
\end{proof}

We now show that the minor iteration sequence terminates at a point
$x^{k,\ell}$ satisfying the acceptance test, provided that $x^k$ is
not a solution.
\begin{theorem} \labtag{th:tr:ft}
  Suppose that $\epstol =0$. 
\begin{itemize}
\item[(i)]   If $x^k \notin \cS$, there is an $\ell \ge 0$ such that
  $x^{k,\ell}$ satisfies \eqnok{tr.accept}. 
\item[(ii)] If $x^k \in \cS$, then either Algorithm TR terminates (and verifies that $x^k \in \cS$), or  
$\cQ(x^k) - m_{k,\ell}(x^{k,\ell}) \downarrow 0$.
\end{itemize}
\end{theorem}
\begin{proof}
  Suppose for the moment that the inner iteration sequence is
  infinite, that is, the test \eqnok{tr.accept} always fails. By
  applying Lemma~\ref{lem:tr:ft} recursively, with any constant
  $\bar{\eta}$ satisfying the properties stated in
  Lemma~\ref{lem:tr:ft}, we can identify a sequence of indices $0 <
  \ell_1 < \ell_2 < \dots$ such that
\beqa
\nonumber
\cQ(x^k) - m_{k,\ell_j}(x^{k,\ell_j}) & \le & 
\bar{\eta} \left[ \cQ(x^k) - m_{k,\ell_{j-1}}(x^{k,\ell_{j-1}}) \right] \\
\nonumber
& \le & 
\bar{\eta}^2  \left[ \cQ(x^k) - m_{k,\ell_{j-2}}(x^{k,\ell_{j-2}}) \right] \\
\nonumber
& \vdots & \\
\labtag{minortozero}
& \le & 
\bar{\eta}^j \left[ \cQ(x^k) - m_{k,0}(x^{k,0}) \right].
\eeqa
When $x^k \notin \cS$, we have from Lemma~\ref{lem:trbounds} that
\[
\Delta_{k,\ell} \ge \min( \Delta_{\rm lo}, E_k/4) 
\defeq \bar{\Delta}_{\rm lo} >0, \;\; \mbox{for all $\ell=0,1,2,\dots$},
\]
so the right-hand side of \eqnok{tr.2a} is strictly positive.  Hence
for $j$ sufficiently large, we have that
\[
\cQ(x^k) - m_{k,\ell_j}(x^{k,\ell_j}) \le 
0.5 \min \left( \bar{\Delta}_{\rm lo}, \| x^k-P(x^k) \|_{\infty} \right)
\frac{\cQ(x^k) - \cQ^*}{\| x^k - P(x^k) \|_{\infty}}.
\]
But this inequality contradicts \eqnok{lem:tr:inequalities}, proving (i).

For the case of $x^k \in \cS$, there are two possibilities. If
the inner iteration sequence terminates finitely at some $x^{k,\ell}$,
we have $\cQ(x^k) - m_{k,\ell}(x^{k,\ell}) = 0$ and indeed that
\[
m_{k,\ell}(x) \ge \cQ(x^k) = \cQ^*, \;\; 
\mbox{for all $x$ with $\| x-x^k \|_{\infty} \le \Delta_{k,\ell}$}.
\]
Because of \eqnok{mkprop.2a}, we have that $\cQ(x) \ge \cQ(x^k)$ for
all $x$ in a neighborhood of $x^k$, implying that $0 \in \partial
\cQ(x^k)$. Therefore, termination under these circumstances yields a
guarantee that $x^k \in \cS$. When the algorithm does not terminate,
it follows from \eqnok{minortozero} that $\cQ(x^k) -
m_{k,\ell}(x^{k,\ell}) \to 0$. By applying Lemma~\ref{lem:mkl}, we
verify our claim (ii) of monotonic convergence.
\end{proof}

We now prove convergence of Algorithm TR to $\cS$.
\begin{theorem} \labtag{th:tr:conv}
  Suppose that $\epstol=0$. The sequence of major
  iterations $\{ x^k \}$ is either finite, terminating at some $x^k
  \in \cS$, or  is infinite, with the property that $\| x^k - P(x^k)
  \|_{\infty} \to 0$.
\end{theorem}
\begin{proof}
  If the claim does not hold, there are two possibilities. The first
  is that the sequence of major iterations terminates finitely at some 
  $x^k \notin \cS$. However, Theorem~\ref{th:tr:ft} ensures, however, that the
  minor iteration sequence will terminate at some new major iteration
  $x^{k+1}$ under these circumstances, so we can rule out this
  possibility. The second possibility is that the sequence $\{x^k\}$
  is infinite but that there is some $\epsilon >0$ and an infinite
  subsequence of indices $\{ k_j \}_{j=1,2,\dots}$ such that
\[
\| x^{k_j} - P(x^{k_j}) \|_{\infty}  \ge \epsilon, \;\; j=0,1,2,\dots.
\]
Since the sequence $\{ \cQ(x^{k_j}) \}_{j=1,2,\dots}$ is infinite,
decreasing, and bounded below, it converges to some value $\bar{\cQ} >
\cQ^*$.  Moreover, since the entire sequence $\{ \cQ(x^k) \}$ is
monotone decreasing, it follows that $\cQ(x^k) > \bar{\cQ}$ and
therefore
\[
\cQ(x^k) - \cQ^* > \bar{\cQ} - \cQ^* > 0, \;\; k=0,1,2,\dots.
\]
Hence, by boundedness of the subgradients (see \eqnok{def.beta}), we
can identify a constant $\bar{\epsilon}>0$ such that
\[
\| x^k - P(x^k) \|_{\infty} \ge \bar{\epsilon}, \;\; k=0,1,2,\dots.
\]
It follows from \eqnok{def:Ek} that
\beq \labtag{Ekbb}
E_k \ge \bar{\epsilon}, \;\; k=0,1,2,\dots.
\eeq

For each major iteration index $k$, let $\ell(k)$ be the minor
iteration index that passes the acceptance test \eqnok{tr.accept}. By combining \eqnok{tr.accept} with Lemma~\ref{lem:tr:1}, we have that
\[
\cQ(x^k) - \cQ(x^{k+1}) \ge \xi \hat{\epsilon} \min 
\left( \Delta_{k, \ell(k)}, \|x^k - P(x^k) \|_{\infty} \right) 
\ge \xi \hat{\epsilon} \min 
\left( \Delta_{k, \ell(k)}, \bar{\epsilon} \right).
\]
Since  $\cQ(x^k) - \cQ(x^{k+1}) \to 0$, we deduce that
\beq \labtag{poo.8}
\lim_{k \to \infty} \Delta_{k, \ell(k)} = 0.
\eeq
By Lemma~\ref{lem:trbounds} and \eqnok{Ekbb}, we have 
\[
 \Delta_{k, \ell(k)} \ge \min (\Delta_{\rm lo}, \bar{\epsilon}/4) >0, \;\;
k=0,1,2,\dots,
\]
which contradicts \eqnok{poo.8}.  We conclude that the second
possibility (an infinite sequence $\{ x^k \}$ not converging to $\cS$)
cannot occur either, so the proof is complete.
\end{proof}

Finally, we show that the algorithm terminates when $\epstol>0$.
\begin{theorem} \labtag{th:fint}
When $\epstol>0$, Algorithm TR terminates finitely.
\end{theorem}
\begin{proof}
  We show first that the algorithm cannot ``get stuck'' at a
  particular $x^k$, generating an infinite sequence of minor
  iterations at $x^k$ without eventually satisfying either
  \eqnok{conv.test} or the acceptance test \eqnok{tr.accept}.  We see
  from the reasoning in the proof of Theorem~\ref{th:tr:ft} together
  with the monotonicity property of Lemma~\ref{lem:mkl} that an
  infinite sequence of minor iterations must satisfy that
\beq \labtag{fint.1}
\cQ(x^k) - m_{k,\ell}(x^{k,\ell}) \downarrow 0.
\eeq
Since the right-hand side of \eqnok{conv.test} is bounded below by
$\epstol$, the test \eqnok{conv.test} must be
satisfied for some $\ell$. Therefore, the minor iteration
sequence cannot be infinite.

Now consider the other possibility of an infinite sequence of major
iterations $\{ x^k \}_{k=1,2,\dots}$. Since we have 
\[
\cQ(x^k) - m_{k,\ell}(x^{k,\ell}) > \epstol
\]
for all $k$ and $\ell$, and since the acceptance test 
\eqnok{tr.accept} is satisfied at all $k$, we have
\[
\cQ(x^k) - \cQ(x^{k+1}) \ge \xi \epstol >0, \;\;
\makebox{for all $k=0,1,2\dots$}.
\]
But this relation is inconsistent with the fact that $\{ \cQ(x^k) \}$
is bounded below (by $\cQ^*$), so this possibility can also be ruled
out, and the proof is complete.
\end{proof}

\subsection{Discussion} \labtag{sec:tr:discussion}

The algorithm can be modified in various ways without
changing its properties greatly.  For instance, we could replace the
step norm bound in \eqnok{trsub.kl} by a scaled bound of the form
\[
\| S (x-x^k) \|_{\infty} \le \Delta_k,
\]
where $S$ is a diagonal positive definite matrix.  After
this modification, \eqnok{master.kl} remains a linear program.  We
could also use a $1$-norm trust region, at the cost of introducing an
additional variable vector $s$ of the same dimension as $x$.
Specifically, we enforce the constraint $\|x-x^k \|_1 \le \Delta_k$ by
enforcing the following linear constraints:
\[
x-x^k \le s, \sgap x^k-x \le s, \sgap e^Ts \le \Delta_k.
\]
Once again, we obtain a linear programming subproblem, albeit one that
involves more variables than \eqnok{master.kl}

If a $2$-norm trust region is used, we can show by comparing the
optimality conditions for the respective problems that the solution of
the subproblem
\[
\min_x \, m_{k,\ell}(x) \;\; \mbox{subject to} \;Ax=b, \; x \ge 0, \; 
\| x-x^k \|_2 \le \Delta_k
\]
is identical to the solution of
\beq \labtag{trsub.2norm}
\min_x \, m_{k,\ell}(x) + \lambda \| x-x^k \|^2 \;\; 
\mbox{subject to} \;Ax=b, \; x \ge 0,
\eeq
for some $\lambda \ge 0$. 
We can transform \eqnok{trsub.2norm} to a quadratic program in the
same fashion as the transformation of \eqnok{trsub.kl} to
\eqnok{master.kl}. The bundle-trust-region approaches described in
Kiwiel~\cite{Kiw90}, Hirart-Urruty and
Lemar\'echal~\cite[Chapter~XV]{HirL93}, and
Ruszczy{\'n}ski~\cite{Rus86,Rus93} also lead to problems of the form
\eqnok{trsub.2norm}. These approaches manipulate the parameter
$\lambda$ rather than adjusting the trust-region radius, more in the
spirit of the Levenberg-Marquardt method for least-squares problems
than of a true trust-region method. Hence, their analysis differs
somewhat from that of the preceding section. Moreover, although
quadratic programming solvers that exploit the special structure of
the quadratic term in \eqnok{trsub.2norm} have been designed and
implemented (see \cite{Rus86}), we believe that the linear programming
subproblem \eqnok{master.kl} is more appealing from a practical point
of view. Improvements in the efficiency and ease of use of linear
programming software have continued to occur at a rapid pace, and
availability of high-quality software has made it much easier to
implement an efficient algorithm based on \eqnok{master.kl} than would
have been the case if the subproblems had the form
\eqnok{trsub.2norm}.

\section{An Asynchronous Bundle-Trust-Region Method}
\labtag{sec:atr}

In this section we present an asynchronous, parallel version of the
trust-region algorithm of the preceding section and analyze its
convergence properties.

\subsection{Algorithm ATR} \labtag{sec:atr:atr}

We now define a variant of the method of Section~\ref{sec:tr} that
allows the partial sums $\cQ_{[j]}, j=1,2,\dots,T$ \eqnok{thetaj} and
their associated cuts to be evaluated simultaneously for different
values of $x$. We generate candidate iterates by solving trust-region
subproblems centered on an ``incumbent'' iterate, which (after a
startup phase) is the point $x^I$ that, roughly speaking, is the best
among those visited by the algorithm whose function value $\cQ(x)$ is
fully known.

By performing evaluations of $\cQ$ at different points concurrently,
we relax the strict synchronicity requirements of Algorithm TR, which
requires $\cQ(x^k)$ to be evaluated fully before the next candidate
$x^{k+1}$ is generated.  The resulting approach, which we call
Algorithm ATR (for ``asynchronous TR''), is more suitable for
implementation on computational grids of the type we consider here.
Besides the obvious increase in parallelism that goes with evaluating
several points at once, there is no longer a risk of the entire
computation being help up by the slow evaluation of one of the partial
sums $\cQ_{[j]}$ on a recalcitrant worker. Algorithm ATR has similar
theoretical properties to Algorithm TR, since the mechanisms for
accepting a point as the new incumbent, adjusting the size of the
trust region, and adding and deleting cuts are all similar to the
corresponding mechanisms in Algorithm TR.

Algorithm ATR maintains a ``basket'' $\cB$ of at most $K$ points for
which the value of $\cQ$ and associated subgradient information is
partially known. When the evaluation of $\cQ(x^q)$ is completed for a
particular point $x^q$ in the basket, it is installed as the new
incumbent if (i) its objective value is smaller than that of the
current incumbent $x^I$; and (ii) it passes a trust-region acceptance
test like \eqnok{tr.accept}, with the incumbent {\em at the time $x^q$
  was generated} playing the role of the previous major iteration in
Algorithm TR.  Whether $x^q$ becomes the incumbent or not, it is
removed from the basket. 

When a vacancy arises in the basket, we may generate a new point by
solving a trust-region subproblem similar to \eqnok{trsub.kl},
centering the trust region at the current incumbent $x^I$.  During the
startup phase, while the basket is being populated, we wait until the
evaluation of some other point in the basket has reached a certain
level of completion (that is, until a proportion $\sigma \in (0,1]$ of
the partial sums \eqnok{thetaj} and their subgradients have been
evaluated) before generating a new point. We use a logical variable
${\tt speceval}_q$ to indicate when the evaluation of $x^q$ passes the
specified threshold and to ensure that $x^q$ does not trigger the
evaluation of more than one new iterate. (Both $\sigma$ and ${\tt
  speceval}_q$ play a similar role in Algorithm ALS.) After the
startup phase is complete (that is, after the basket has been filled),
vacancies arise only after evaluation of an iterate $x^q$ is
completed.

We use $m(\cdot)$ (without
subscripts) to denote the model function for $\cQ(\cdot)$. When
generating a new iterate, we use whatever cuts are stored at the
time to define $m$.  When solved around the incumbent $x^I$
with trust-region radius $\Delta$, the subproblem is as follows: 
\beq
\labtag{trsub.atr1} \mbox{\tt trsub$(x^I, \Delta)$:} \;\; \min_x \,
m(x) \;\; \mbox{subject to} \;Ax=b, \; x \ge 0, \; \| x- x^I
\|_{\infty} \le \Delta.  
\eeq 
We refer to $x^I$ as the {\em parent incumbent} of the solution of
\eqnok{trsub.atr1}.

In the following description, we use $k$ to index the successive
points $x^k$ that are explored by the algorithm, $I$ to denote the
index of the incumbent, and $\cB$ to denote the basket.  We use $t_k$
to count the number of partial sums $\cQ_{[j]}(x^k)$, $j=1,2,\dots,T$
that have been evaluated so far.

Given a starting guess $x^0$, we initialize the algorithm by setting
the dummy point $x^{-1}$ to $x^0$, setting the incumbent index $I$ to
$-1$, and setting the initial incumbent value $\cQ^I =\cQ^{-1}$ to
$\infty$. The iterate at which the first evaluation is completed
becomes the first ``serious'' incumbent.

We now outline some other notation used in specifying Algorithm ATR:
\bi

\item[$\cQ^I$:] The objective value of the incumbent $x^I$, except in
the case of $I=-1$, in which case $\cQ^{-1} = \infty$.

\item[$I_q$:] The index of the parent incumbent of $x^q$, that is, the
  incumbent index $I$ at the time that $x^q$ was generated from
  \eqnok{trsub.atr1}. Hence, $\cQ^{I_q} = \cQ(x^{I_q})$ (except when
  $I_q=-1$; see previous item).

\item[$\Delta_q$:] The value of the trust-region radius $\Delta$ used 
when solving for  $x^q$.

\item[$\Delta_{\rm curr}$:] Current value of the trust-region
radius. When it comes time to solve \eqnok{trsub.atr1} to obtain a new
iterate $x^q$, we set $\Delta_q \leftarrow \Delta_{\rm curr}$.

\item[$m^q$:] The optimal value of the objective function $m$ in the
subproblem {\tt trsub}$(x^{I_q}, \Delta_q)$ \eqnok{trsub.atr1}.

\ei

Our strategy for maintaining the model closely follows that of
Algorithm TR. Whenever the incumbent changes, we have a fairly free
hand in deleting the cuts that define $m$, just as we do after
accepting a new major iterate in Algorithm TR. If the incumbent does
not change for a long sequence of iterations (corresponding to a long
sequence of minor iterations in Algorithm TR), we can still delete
``stale'' cuts that represent information in $m$ that has likely been
superseded (as quantified by a parameter $\eta \in [0,1)$). The
following version of Procedure Model-Update, which applies to
Algorithm ATR, takes as an argument the index $k$ of the latest
iterate generated by the algorithm. It is called after the evaluation
of $\cQ$ at an earlier iterate $x^q$ has just been completed, but
$x^q$ does {\em not} meet the conditions needed to become the new
incumbent.
\btab
\> {\bf Procedure Model-Update} $(k)$ \\
\> {\bf for each} optimality cut defining $m$\\
\>\> {\tt possible\_delete}  $\leftarrow$ {\tt true}; \\
\>\> {\bf if} the cut was generated at the parent incumbent $I_k$ of $k$\\
\>\>\> {\tt possible\_delete}  $\leftarrow$ {\tt false}; \\
\>\> {\bf else if} the cut was active at the solution $x^k$ of 
{\tt trsub}$(x^{I_k},\Delta_k)$ \\
\>\>\> {\tt possible\_delete}  $\leftarrow$ {\tt false}; \\
\>\> {\bf else if} the cut was generated at an earlier 
iteration $\bar{\ell}$ \\
\>\>\>\> such that $I_{\bar{\ell}} = I_k \neq -1$ and
\etab
\beq \labtag{atr.cut.delete.criterion}
\cQ^{I_k} - m^k > \eta [ \cQ^{I_k} - m^{\bar{\ell}} ]
\eeq
\btab
\>\>\> {\tt possible\_delete}  $\leftarrow$ {\tt false}; \\
\>\> {\bf end (if)} \\
\>\> {\bf if} {\tt possible\_delete} \\
\>\>\> possibly delete the cut; \\
\> {\bf end (for each)}
\etab

Our strategy for adjusting the trust region $\Delta_{\rm curr}$
also follows that of Algorithm TR. The differences arise from the fact
that between the time an iterate $x^q$ is generated and its function
value $\cQ(x^q)$ becomes known, other adjustments of $\Delta_{\rm
current}$ may have occurred, as the evaluation of intervening iterates
is completed. The version of Procedure Reduce-$\Delta$  for 
Algorithm ATR is as follows.
\btab
\> {\bf Procedure Reduce-$\Delta(q)$} \\
\> {\bf if} $I_q = -1$ \\
\>\> return; \\
\> evaluate
\etab
\beq \labtag{atr.reduce.delta.2}
\rho = {\min(1,\Delta_q)} 
\frac{\cQ(x^q) - \cQ^{I_q}}{\cQ^{I_q}  - m^q};
\eeq
\btab
\> {\bf if} $\rho>0$ \\
\>\> {\tt counter} $\leftarrow$ {\tt counter}$+1$; \\
\> {\bf if} $\rho>3$ {\bf or} 
({\tt counter} $\ge 3$ {\bf and} $\rho \in (1,3]$) \\
\>\> set $\Delta_q^+ \leftarrow \Delta_q / \min(\rho,4)$; \\
\>\> set 
$\Delta_{\rm curr} \leftarrow \min(\Delta_{\rm curr}, \Delta_q^+)$; \\
\>\> reset {\tt counter} $\leftarrow 0$; \\
\> return.
\etab

The protocol for increasing the trust region after a successful step
is based on \eqnok{tr.incr.1}, \eqnok{tr.incr.3}. If on completion of
evaluation of $\cQ(x^q)$, the iterate $x^q$ becomes the new incumbent,
then we test the following condition:
\beq \labtag{atr.incr.1}
\cQ(x^q) \le \cQ^{I_q} - 0.5 (\cQ^{I_q} - m^q) \;\; \mbox{and} \;\;
\| x^q - x^{I_q} \|_{\infty} = \Delta_q.
\eeq
If this condition is satisfied, we set
\beq \labtag{atr.incr.3}
\Delta_{\rm curr} \leftarrow \max(\Delta_{\rm curr},
 \min (\Delta_{\rm hi}, 2 \Delta_q) ).
\eeq

The convergence test is also similar to the test \eqnok{conv.test}
used for Algorithm TR. We terminate if, on generation of a new iterate
$x^k$, we find that
\beq \labtag{conv.test.atr}
\cQ^I - m^k \le \epstol (1+|\cQ^I|).
\eeq

We now specify the four key routines of the Algorithm ATR, which serve
a similar function to the four main routines of Algorithm ALS. As in
the earlier case, we assume for simplicity of description that each
task consists of evaluation of the function and a subgradient for
a single cluster (although in practice we may bundle more than one
cluster into a single task). The routine {\tt partial\_evaluate}
executes on worker processors, while the other three routines execute
on the master processor.

\btab
\>{\bf ATR:} \ \ {\tt  partial\_evaluate}$(x^q,q,j,\cQ_{[j]}(x^q),g_j)$ \\
\> Given $x^q$, index  $q$, and  partition number $j$, 
evaluate $\cQ_{[j]}(x^q)$ from \eqnok{thetaj} \\
\>\> together with a partial subgradient $g_j$ from \eqnok{subg.Qj}; \\
\> Activate {\tt act\_on\_completed\_task}$(x^q,q,j,\cQ_{[j]}(x^q),g_j)$
on the master processor.
\etab

\medskip

\btab
\> {\bf ATR:} \ \ {\tt  evaluate}$(x^q,q)$ \\
\> {\bf for} $j=1,2,\dots, T$ (possibly concurrently) \\
\>\> {\tt partial\_evaluate}$(x^q,q,j,\cQ_{[j]}(x^q), g_j)$; \\
\> {\bf end (for)}
\etab

\medskip

\btab
\> {\bf ATR:} \ \ {\tt initialization}$(x^0)$ \\
\> choose $\xi \in (0,1/2)$, trust region upper bound
$\Delta_{\rm hi}>0$; \\
\> choose synchronicity parameter $\sigma \in (0,1]$; \\
\> choose maximum basket size $K>0$; \\
\> choose $\Delta_{\rm curr} \in (0, \Delta_{\rm hi}]$, 
{\tt counter} $\leftarrow 0$; $\cB \leftarrow \emptyset$; \\
\> $I \leftarrow -1$; $x^{-1} \leftarrow x^0$; $\cQ^{-1} \leftarrow \infty$; 
$I_0 \leftarrow -1$; \\
\> $k \leftarrow 0$; 
${\tt speceval}_0 \leftarrow {\tt  false}$; 
$t_0 \leftarrow 0$; \\
\> {\tt evaluate}$(x^0,0)$.
\etab

\medskip

\btab
\> {\bf ATR:} \ \  
{\tt act\_on\_completed\_task}$(x^q,q,j,\cQ_{[j]}(x^q),g_j))$ \\
\> $t_q \leftarrow t_q+1$; \\
\> add $\cQ_{[j]}(x^q)$ and cut $g_j$ to the model $m$; \\
\> {\tt basketFill} $\leftarrow$ {\tt  false}; 
{\tt basketUpdate} $\leftarrow$ {\tt  false}; \\
\> {\bf if} $t_q=T$ (* evaluation of $\cQ(x^q)$ is complete *) \\
\>\> {\bf if} $\cQ(x^q) < \cQ^I$ and (${I_q}=-1$ or
$\cQ(x^q) \le \cQ^{I_q} - \xi (\cQ^{I_q} - m^q)$) \\
\>\>\> (* make $x^q$ the new incumbent *) \\
\>\>\>  $I \leftarrow q$;  $\cQ^I \leftarrow \cQ(x^I)$; \\
\>\>\> possibly increase  $\Delta_{\rm curr}$ according to 
\eqnok{atr.incr.1} and \eqnok{atr.incr.3}; \\
\>\>\> modify the model function by possibly deleting cuts not arising \\
\>\>\>\> from the evaluation of $\cQ(x^q)$; \\
\>\> {\bf else} \\
\>\>\> call Model-Update$(k)$; \\
\>\>\> call Reduce-$\Delta(q)$ to update $\Delta_{\rm curr}$; \\
\>\> {\bf end (if)} \\
\>\> $\cB \leftarrow \cB \backslash \{ q \}$; \\
\>\> {\tt basketUpdate} $\leftarrow$ {\tt true}; \\

\> {\bf else if } 
 $t_q \ge \sigma T$ {\bf and} $| \cB| <K$ {\bf and} not ${\tt speceval}_q$ \\
\>\> (* basket-filling phase: enough partial sums have been evaluated at $x^q$ 
 \\
\>\>\> to trigger calculation of a new candidate iterate *) \\
\>\> ${\tt speceval}_q \leftarrow ${\tt true}; 
{\tt basketFill} $\leftarrow$ {\tt true}; \\
\> {\bf end (if)} \\

\> {\bf if } {\tt basketFill} {or}
{\tt basketUpdate} \\
\>\> $k \leftarrow k+1$; 
set $\Delta_k \leftarrow \Delta_{\rm curr}$; set $I_k \leftarrow I$; \\
\>\> solve {\tt trsub}$(x^I,\Delta_k)$ to obtain $x^k$; \\
\>\> $m^k \leftarrow m(x^k)$; \\
\>\> {\bf if} \eqnok{conv.test.atr} holds \\
\>\>\> STOP;  \\
\>\> $\cB \leftarrow \cB \cup \{ k \}$; \\
\>\> ${\tt speceval}_k \leftarrow${\tt false}; $t_k \leftarrow 0$; \\
\>\> {\tt evaluate}$(x^k,k)$; \\
\> {\bf end (if)}

\etab

It is not generally true that the first $K$ iterates $x^0, x^1, \dots,
x^{K-1}$ generated by the algorithm are all basket-filling
iterates. Often, an evaluation of some iterate is completed before the
basket has filled completely, so a ``basket-update'' iterate is used
to generate a replacement for this point. Since each basket-update
iterate does not change the size of the basket, however, the number of
basket-filling iterates that are generated in the course of the
algorithm is exactly $K$.

\subsection{Analysis of Algorithm ATR} \labtag{sec:atr:analysis}

We now analyze Algorithm ATR, showing that its convergence properties
are similar to those of Algorithm TR. Throughout, we make the
following assumption:
\beq \labtag{all.tasks.completed}
\mbox{Every task is completed after a finite  time}. 
\eeq

The analysis follows closely that of Algorithm TR presented in
Section~\ref{sec:tr:analysis}. We state the analogues of all the
lemmas and theorems from the earlier section, incorporating the
changes and redefinitions needed to handle Algorithm ATR. Most of the
details of the proofs are omitted, however, since they are similar to
those of the earlier results.

We start by defining the level set within which the points and
incumbents generated by ATR lie.
\begin{lemma} \labtag{lem:atr1.1}
All incumbents $x^I$ generated by ATR lie in $\cL(\cQ_{\rm max})$,
whereas all points $x^k$ considered by the algorithm lie in
$\cL(\cQ_{\rm max}; \Delta_{\rm hi})$, where $\cL(\cdot)$ and
$\cL(\cdot;\cdot)$ are defined by \eqnok{def.ls} and \eqnok{def.lsn},
respectively, and $\cQ_{\rm max}$ is defined by
\[
\cQ_{\rm max} \defeq \sup \{ \cQ(x) \, | \, 
\| x-x^0 \| \le \Delta_{\rm hi} \}.
\]
\end{lemma}
\begin{proof}
  Consider first what happens in ATR before the  first function
  evaluation is complete. Up to this point, all the iterates $x^k$ in
  the basket are generated in the basket-filling part and therefore
  satisfy $\| x^k-x^0 \| \le \Delta_k \le \Delta_{\rm hi}$, with
  $\cQ^{I_k} = \cQ^{-1} = \infty$.
  
  When the first evaluation is completed (by $x^k$, say), it trivially
  passes the test to be accepted as the new incumbent. Hence, the
  first noninfinite incumbent value becomes $\cQ^I = \cQ(x^k)$, and
  by definition we have $\cQ^I \le \cQ_{\rm max}$.  Since all later
  incumbents must have objective values smaller than this first
  $\cQ^I$, they all must lie in the level set $\cL(\cQ_{\rm max})$,
  proving our first statement.

All points $x^k$ generated within {\tt act\_on\_completed\_task} lie
within a distance $\Delta_k \le \Delta_{\rm hi}$ either of $x^0$ or of
one of the later incumbents $x^I$. Since all the incumbents, including
$x^0$, lie in $\cL(\cQ_{\rm max})$, we conclude that the second claim
in the theorem is also true.

\end{proof}

Analogously with $\beta$ \eqnok{def.beta}, we define a bound on the
subgradients over the set $\cL(\cQ_{\rm max}; \Delta_{\rm hi})$ as
follows:
\beq \labtag{def.barbeta}
\bar{\beta} = \sup \{ \| g \|_1 \, | \, g \in \partial \cQ(x), \, 
\mbox{for some $x \in \cL(\cQ_{\rm max};\Delta_{\rm hi})$} \}.
\eeq

The next result is analogous to Lemma~\ref{lem:mkl}. It shows that for
any sequence of iterates $x^k$ for which the parent incumbent $x^I_k$
is the same, the optimal objective value in {\tt trsub}$(x^{I_k},
\Delta_k)$ is monotonically increasing.
\begin{lemma} \labtag{lem:mkl.atr}
Consider any contiguous subsequence of iterates $x^{k}$,
$k=k_1,k_1+1,\dots, k_2$ for which the parent incumbent is identical;
that is, $I_{k_1}=I_{k_1+1}= \cdots = I_{k_2}$. Then we have
\[
m^{k_1} \le m^{k_1+1} \le \cdots \le m^{k_2}.
\]
\end{lemma}
\begin{proof}
We select any $k=k_1, k_1+1, \dots, k_2-1$ and 
prove that $m^k \le m^{k+1}$.
Since $x^k$ and $x^{k+1}$ have the same parent incumbent ($x^I$, say),
no new incumbent has been accepted between the generation of these two
iterates, so the wholesale cut deletion that may occur with the
adoption of a new incumbent cannot have occurred.  There may, however,
have been a call to {\tt Model-Update}$(k)$. The
first ``else if'' clause in {\tt Model-Update} would have ensured that
cuts active at the solution of {\tt trsub}$(x^I, \Delta_k)$ were still
present in the model when we solved {\tt trsub}$(x^I, \Delta_{k+1})$ to
obtain $x^{k+1}$. Moreover, since no new incumbent was accepted,
$\Delta_{\rm curr}$ cannot have been increased, and we have
$\Delta_{k+1} \le \Delta_k$. We now use the same argument as in the
proof of Lemma~\ref{lem:mkl} to deduce that $m^{k} \le m^{k+1}$.
\end{proof}

The following result is analogous to Lemma~\ref{lem:tr:1}. We omit the
proof, which modulo the change in notation is identical to the earlier
result.
\begin{lemma} \labtag{lem:atr:1}
For all $k=0,1,2,\ldots$ such that $I_k \neq -1$,  we have that 
\begin{subequations} \labtag{lem:atr:inequalities}
\beqa 
\labtag{atr.2a}
\cQ^{I_k} - m^k & \ge  &
\min \left( \Delta_{k}, \| x^{I_k} - P(x^{I_k})\|_{\infty} \right)
\frac{\cQ^{I_k} - \cQ^*}{\| x^{I_k} - P(x^{I_k}) \|_{\infty}} \\
\labtag{atr.2b}
& \ge &
\hat{\epsilon} \min 
\left( \Delta_{k}, \| x^{I_k} - P(x^{I_k})\|_{\infty} \right),
\eeqa
\end{subequations}
where $\hat{\epsilon}>0$ is defined in \eqnok{weak.sharp}.
\end{lemma}

The following analogue of Lemma~\ref{lem:trbounds} requires a slight
redefinition of the quantity $E_k$ from \eqnok{def:Ek}. We now
define it to be the closest approach by an {\em incumbent} to the
solution set, up to and including iteration $k$; that is, 
\beq \labtag{def:Ek:atr}
E_k \defeq \min_{\bar{k} = 0,1,\dots, k; I_{\bar{k}} \neq -1} 
\| x^{I_{\bar{k}}} - P(x^{I_{\bar{k}}}) \|_{\infty}.
\eeq
We also omit the proof of the following result, which, allowing for
the change of notation, is almost identical to that of
Lemma~\ref{lem:trbounds}.
\begin{lemma} \labtag{lem:trbounds:atr}
  There is a constant $\Delta_{\rm lo} >0$ such that for all trust
  regions $\Delta_{k}$ used in the course of Algorithm ATR, we
  have
\[
\Delta_{k} \ge \min( \Delta_{\rm lo}, E_k/4).
\]
\end{lemma}
The value of $\Delta_{\rm lo}$ that works in this case is $\Delta_{\rm
  lo} = (1/4) \min(1, \hat{\epsilon}/\bar{\beta}, \Delta_{\rm hi})$,
where $\bar{\beta}$ comes from \eqnok{def.barbeta}.

There is also an analogue of Lemma~\ref{lem:tr:ft} that shows that if
the incumbent remains the same for a number of consecutive iterations,
the gap between incumbent objective value and model function decreases
significantly as the iterations proceed.
\begin{lemma} \labtag{lem:atr:ft}
  Let $\epstol=0$ in Algorithm ATR, and let $\bar{\eta}$ be
  any constant satisfying $0<\bar{\eta}<1$, $\bar{\eta}>\xi$,
  $\bar{\eta} \ge \eta$. Choosing any index $k_1$ with $I_{k_1} \neq
  -1$, we have either that the incumbent $I_{k_1}=I$ is eventually
  replaced by a new incumbent or that there is an iteration
  $k_2>k_1$ such that 
\beq \labtag{atr.6}
\cQ^{I} - m^{k_2} \le \bar{\eta} \left[
\cQ^{I} - m^{k_1} \right].
\eeq
\end{lemma}
The proof of this result follows closely that of its antecedent
Lemma~\ref{lem:tr:ft}. The key is in the construction of the
Model-Update procedure. As long as 
\beq \labtag{atr.7}
\cQ^I - m^k > \eta [\cQ^I - m^{k_1}], \;\; \mbox{for $k \ge k_1$, where
$I=I_{k_1} = I_k$},
\eeq
none of the cuts generated during the evaluation of $\cQ(x^q)$ for any
$q=k_1, k_1+1, \dots, k$ can be deleted. The proof technique of
Lemma~\ref{lem:tr:ft} can then be used to show that the successive
iterates $x^{k_1}, x^{k_1+1}, \dots$ cannot be too closely spaced if
the condition \eqnok{atr.7} is to hold and if all of them fail to
satisfy the test to become a new incumbent. Since they all belong
to a box of finite size centered on $x^I$, there can be only finitely
many of these iterates. Hence, either a new incumbent is adopted
at some iteration $k \ge k_1$ or  condition \eqnok{atr.6} is 
eventually satisfied.

We now show that the algorithm cannot ``get stuck'' at a nonoptimal
incumbent. The following result is analogous to
Theorem~\ref{th:tr:ft}, and its proof relies on the earlier results in
exactly the same way.
\begin{theorem} \labtag{th:atr:ft}
  Suppose that $\epstol =0$. 
\begin{itemize}
\item[(i)] If $x^I \notin \cS$, then this incumbent is replaced by a
new incumbent after a finite time.
\item[(ii)] If $x^I \in \cS$, then either Algorithm ATR terminates
(and verifies that $x^I \in \cS$), or $\cQ^I - m^k \downarrow 0$
as $k \to \infty$.
\end{itemize}
\end{theorem}

We conclude with the result that shows convergence of the sequence of
incumbents to $\cS$. Once again, the logic of proof follows that of
the synchronous analogue Theorem~\ref{th:tr:conv}.
\begin{theorem} \labtag{th:atr:conv}
  Suppose that $\epstol=0$. The sequence of incumbents
  $\{ x^{I_k} \}_{k=0,1,2,\dots}$ is either finite, 
terminating at some $x^I \in \cS$ or is infinite with 
 the property that $\| x^{I_k} - P(x^{I_k})
  \|_{\infty} \to 0$.
\end{theorem}

\section{Implementation on Computational Grids} \labtag{sec:grids}

We now describe some salient properties of the computational
environment in which we implemented the algorithms, namely, a
computational grid running the Condor system and the MW runtime
support library.

\subsection{Properties of Grids} \labtag{sec:grids:intro}

The term ``grid computing'' (synonymously ``metacomputing'') is
generally used to describe parallel computations on a geographically
distributed, heterogeneous computing platform. Within this framework
there are several variants of the concept. The one of interest here is
a parallel platform made up of shared workstations, nodes of PC
clusters, and supercomputers.  Although such platforms are potentially
powerful and inexpensive, they are difficult to harness for productive
use, for the following reasons:
\bi
\item Poor communications properties. Latencies between the processors
  may be high, variable, and unpredictable.
  
\item Unreliability. Resources may disappear without notice. A
  workstation performing part of our computation may be reclaimed by
  its owner and our job terminated.
  
\item Dynamic availability. The pool of available processors grows and
shrinks during the computation, according to the claims of other users
and scheduling considerations at some of the nodes.

\item Heterogeneity. Resources may vary in their operational
characteristics (memory, swap space, processor speed, operating
system).

\ei
In all these respects, our target platform differs from conventional
multiprocessor platforms (such as IBM SP or SGI Origin machines) and
from Linux clusters.

\subsection{Condor} \labtag{sec:grids:condor}

Our particular interest is in grid computing platforms based on the
Condor system~\cite{condor}, which manages distributively owned
collections (``pools'') of processors of different types, including
workstations, nodes from PC clusters, and nodes from conventional
multiprocessor platforms. When a user submits a job, the Condor system
discovers a suitable processor for the job in the pool, transfers the
executable and starts the
job on that processor. It traps system calls (such as input/output
operations), referring them back to the submitting workstation,
and checkpoints the state of the job periodically. It also migrates the
job to a different processor in the pool if the current host becomes
unavailable for any reason (for example, if the workstation is
reclaimed by its owner).  Condor managed
processes can communicate through a Condor-enabled version of PVM
\cite{PVMbook} or by using Condor's I/O trapping to write into and
read from a series of shared files.

\subsection{Implementation in MW} \labtag{sec:grids:mw}

MW (see Goux, Linderoth, and Yoder~\cite{GouLY00} and Goux et
al.~\cite{GouKLY00}) is a runtime support library that facilitates
implementation of parallel master-worker applications on computational
grids. To implement MW on a particular computational grid, a grid
programmer must reimplement a small number of functions to perform
basic operations for communications between processors and management
of computational resources. These functions are encapsulated in the
MWRMComm class. Of more relevance to the current paper is the other
side of MW, the application programming interface presented to the
application programmer. This interface takes the form of a set of
three C$++$ abstract classes that must be reimplemented in a way that
describes the particular application. These classes, named MWDriver,
MWTask, and MWWorker, contain a total of ten methods for which the
user must supply implementations. We describe these methods briefly,
indicating how they are implemented for the particular case of the ATR
and ALS algorithms.

\paragraph{MWDriver.}

This class is made up of methods that execute on the submitting
workstation, which acts as the master processor. It contains the
following four C$++$ pure virtual functions. (Naturally, other methods
can be defined as needed to implement parts of the algorithm.)
\begin{itemize}
\item {\tt get\_userinfo}: Processes command-line arguments and does
  basic setup. In our applications this function reads a command file
  to set various parameters, including convergence tolerances, number
  of scenarios, number of partial sums to be evaluated in each task,
  maximum number of worker processors to be requested, initial trust
  region radius, and so on. It calls the routines that read and store
  the problem data files, and it reads the initial point, if one is
  supplied.  It also performs the operations specified in the {\tt
    initialization} routine of Algorithms ALS and ATR, except for the
  final {\tt evaluate} operation, which is handled by the next
  function.

\item {\tt setup\_initial\_tasks}: Defines the initial pool of tasks.
  In the case of Algorithms ALS and ATR, this function corresponds to
  a call to {\tt evaluate} at $x^0$.

\item {\tt pack\_worker\_init\_data}: Packs the initial data to be
  sent to each worker processor when it joins the pool. In our case,
  the information contained in the input files for the stochastic
  programming problem is sent to each worker.  When the worker
  subsequently receives a task requiring it to solve a number of
  second-stage scenarios, it can use the original input data to
  generate the particular data for its assigned set of scenarios.
  By loading each new worker with the problem data, we avoid having to
  subsequently pass a complete set of data for every scenario in every
  task.

\item {\tt act\_on\_completed\_task}: Is called every time
  a task finishes, to process the results of the task and to take any
  actions arising from these results.  See Algorithms ALS and ATR for
  our definition of this function in our applications. 

\end{itemize}

The MWDriver base class performs many other operations associated with
handling worker processes that join and leave the computation,
assigning tasks to appropriate workers, rescheduling tasks when their
host workers disappear without warning, and keeping track of
performance data for the run. All this complexity is hidden from the
application programmer.

\paragraph{MWTask.}

The MWTask is the abstraction of a single task. It holds both the data
describing that task and the results obtained by executing the task.
The user must implement four functions for packing and unpacking this
data and results between master and workers into simple data
structures that can be communicated between master and workers using
the appropriate primitives for the particular computational grid
platform on which MW is implemented. In most of the results reported
in Section~\ref{sec:results}, the message-passing facilities of
Condor-PVM were used to perform the communication.  By simply changing
compiler directives, the same algorithmic code can also be implemented
on an alternative communication protocol that uses shared files to
pass messages between master and workers. The large run reported in
the next section used this version of the code.

In our applications, each task evaluates the partial sum
$\cQ_{[j]}(x)$ and a subgradient for a given number of clusters. The
task is described by a range of scenario indices for each cluster in
the task and by a value of the first-stage variables $x$. The results
consist of the function and subgradient for each of the clusters
in the task.

\paragraph{MWWorker.}

The MWWorker class is the core of the executable that runs on each
worker. The user must implement two pure virtual functions:

\begin{itemize}
\item {\tt unpack\_init\_data}: Unpacks the initial information passed
  to the worker by the MWDriver function {\tt
    pack\_worker\_init\_data()} when the worker joins the pool. (See
  the discussion of {\tt pack\_worker\_init\_data} in the MWDriver class.)

\item {\tt execute\_task}: Executes a single task.
\end{itemize}

After initializing itself, using the information passed to it by the
master, the worker process sits in a loop, waiting for tasks to be
sent to it. When it detects a new task, it calls {\tt execute\_task}
to compute the results. It passes the results back to the worker by
using the appropriate function from the MWTask class, and then returns
to its wait loop. The wait loop terminates when the master sends a
termination message. In our applications, the {\tt execute\_task()}
function formulates the second-stage linear programs in its clusters
by using the information in the task definition and the data passed to
the worker on initialization. It then calls the linear programming
solvers SOPLEX or CPLEX
 to solve these linear programs, and
uses the dual solutions to calculate the subgradient for each cluster.

\section{Computational Results} \labtag{sec:results}

We now report on computational experiments obtained with
implementations of the ALS, TR, and ATR algorithms using MW on the
Condor system. After describing some further details of the
implementations and the experiments, we discuss our choices for the
various algorithmic parameters and how these were varied between runs.
We then tabulate and discuss the results.

\subsection{Implementations and Experiments}
\label{sec:results:details}

As noted earlier, we used the Condor-PVM implementation of MW for most
of the the runs reported here.
Most of the computational time is taken up with solving linear
programming problems, both by the master process (in the {\tt
  act\_on\_completed\_task} function) and in the tasks, which solve
clusters of second-stage linear programs. We used the CPLEX simplex
solver on the master processor and the SOPLEX public-domain simplex
code (see Wunderling~\cite{soplex}) on the workers. SOPLEX is somewhat
slower in general, but since most of the machines in the Condor pool
do not have CPLEX licenses, there was little alternative but
to use a public-domain code.

We ran most of our experiments on the Condor pool at the University of
Wisconsin, sometimes using Condor's flocking mechanism to augment this
pool with processors from other sites. The other sites included the
University of New Mexico, Columbia University, and the Linux cluster
Chiba City at Argonne National Laboratory. The architectures included
PCs running Linux, and PCs and Sun workstations running different
versions of Solaris. The number of workers available for our use
varied dramatically between and during each set of trials, because of
the differing priorities of the two accounts we used, the variation of
our priority during each run, the number and priorities of other users
of the Condor pool at the time, and the varying number of machines
available to the pool.  The latter number tends to be larger during
the night, when owners of the individual workstations are less likely
to be using them.  The master process was run on a Linux machine in
some experiments and an Intel Solaris machine in other cases.

The input files for the problems reported here were in SMPS format
(see Birge et al.~\cite{BirDGGKW87} and Gassmann and
Schweitzer~\cite{GasS97}). We considered two-stage stochastic linear
programs in which the number of scenarios is finite but extremely
large. We used Monte Carlo sampling to obtain approximate problems
with a specified number $N$ of second-stage scenarios. Brief
descriptions of the test problems can be found at \cite{Hol97}.
In each experiment, we supplied a starting point to the code, obtained
from the solution of a different sampled instance of the same problem.
The function value of the starting point was therefore quite close to
the optimal objective value.

\subsection{Critical Parameters}
\label{sec:results:parameters}

As part of the initialization procedure (implemented by the {\tt
  get\_userinfo} function in the MWDriver class), the code reads an
input file in which various parameters are specified. Several
parameters, such as those associated with modifying the size of the
trust region, have fixed values that we have discussed already in the
text. Others are assigned the same values for all algorithms and all
experiments, namely,
\[
\epsilon_{\rm tol} = 10^{-5}, \sgap
\Delta_{\rm hi} = 10^3, \sgap
\Delta_{0,0} = \Delta_0 = 1, \sgap
\xi = 10^{-4}.
\]
We also set $\eta= 0$ in the Model-Update functions in both TR and
ATR. In TR, this choice has the effect of not allowing deletion of
cuts generated during any major iterations, until a new major iterate
is accepted. In ATR, the effect is to not allow deletion of cuts that
are generated at points whose parent incumbent is still the incumbent.
Even among cuts for which {\tt possible\_delete} is still true at the
final conditional statement of the Model-Update procedures, we do not
actually delete the cuts until they have been inactive at the solution
of the trust-region subproblem for a specified number of consecutive
iterations. For TR, we delete the cut if it has been inactive for more
than 100 consecutive minor iterations, while in ATR we delete the cut
if it was last active at subproblem $\ell$, where $\ell < k-100$ and
$k$ is the current iteration index. Our cut deletion strategy is
therefore not at all parsimonious; it tends to lead to subproblems
\eqnok{trsub.kl} and \eqnok{trsub.atr1} with fairly large numbers of
cuts. In most cases, however, the storage required for these cuts and
the time required to solve the subproblems remain reasonable. We
discuss the exceptions below.

The synchronicity parameter $\sigma$, which arises in Algorithms ALS
and ATR and which specifies the proportion of clusters from a
particular point that must be evaluated in order to trigger evaluation
of a new candidate solution, is varied between $.5$ and $1.0$ in our
experiments.  The size $K$ of the basket $\cB$ is varied between $1$
and $14$. For each problem, the number $T$ of clusters is also varied
in a manner described in the tables, as is the number of tasks into
which the second-stage calculations are divided, which we denote by
$C$. Note that the number of second-stage LPs per chunk is therefore
$N/C$ while the number per cluster is $N/T$.

The MW library allows us to specify an upper bound on the number of
workers we request from the Condor pool, so that we can avoid claiming
more workers than we can utilize effectively. We calculate a rough
estimate of this number based on the number of tasks $C$ per
evaluation of $\cQ(x)$ and the basket size $K$.  For instance, the
synchronous TR and LS algorithms can never use more than $C$ worker
processors, since they evaluate $\cQ$ at just one $x$ at a time. In
the case of TR and ATR, we request $\mbox{mid} (25, 200, \lfloor
(K+1)C/2 \rfloor)$
workers.  For ALS, we request $\mbox{mid}(25,200,2C)$ workers.

We have a single code that implements all four algorithms LS, ALS, TR,
and ATR, using logical branches within the code to distinguish between
the L-shaped and trust-region variants.  There is no distinction in
the code between the two synchronous variants and their asynchronous
counterparts. Instead, by setting $\sigma=1.0$, we force synchronicity
by ensuring that the algorithm considers only one value of $x$ at a
time.

Whenever a worker processor joins the computation, MW sends it a
benchmark task that typifies the type of task it will receive during
the run. In our case, we define the benchmark task to be the solution
of $N/C$ second-stage LPs. The time required for the processor to
solve this task is logged, and we set the ordering policy so as to
ensure that when more than one worker is available to process a
particular task, the task is sent to the worker that logged the
fastest time on the benchmark task.

\subsection{Results: Varying Parameter Choices} \label{sec:results:numbers}

In this section we describe a series of experiments on the same
problem, using different parameter settings, and run under different
conditions on the Condor pool. For these trials, we use the problem
SSN, which arises from a network design application described by Sen,
Doverspike, and Cosares~\cite{SenDC94}. This problem is based on a
graph with 89 arcs, each representing a telecommunications link
between two cities. The first-stage variables represent the
(nonnegative) extra capacity to be added to each of these 89 arcs to
meet an uncertain demand pattern. There is a constraint on the total
added capacity. The demands consist of requests for service between
pairs of nodes in the graph. For each set of requests, a route through
the network of sufficient capacity to meet the requests must be found,
otherwise a penalty term for each request that cannot be satisfied is
added to the objective. The second-stage problems are network flow
problems for calculating the routing for a given set of demand flows.
Each such problem is nontrivial: 706 variables, 175 constraints, and
2284 nonzeros in the constraint matrix. The uncertainty lies in the
fact that the demand for service on each of the 86 pairs is not known
exactly. Rather, there are three to seven possible scenarios for
these demands, all independent of each other, giving a total of about
$10^{70}$ possible scenarios. We use Monte Carlo sampling to obtain a
sampled approximation with $N=10,000$ scenarios. The deterministic
equivalent for this sampled approximation has approximately $1.75
\times 10^6$ constraints and $7.06 \times 10^6$ variables. In all the
runs, we used as starting point the computed solution for a different
sampled approximation---one with $20,000$ scenarios and a different
random seed. The starting point had a function value of approximately
$9.868860$, whereas the optimal objective was approximately
$9.832544$.

In the tables below we list the following information.
\begin{itemize}
\item {\bf points evaluated}. The number of distinct values of the
first-stage variables $x$ generated by solving the master
subproblem---the problem \eqnok{als.subprob} for Algorithm ALS,
\eqnok{trsub.kl} for Algorithm TR, and \eqnok{trsub.atr1} for
Algorithm ATR. 

\item {\bf $| \cB |$}. Maximum size of the basket, also denoted above by $K$.

\item {\bf number of tasks (chunks)}. Denoted above by $C$.

\item {\bf number of clusters}. Denoted above by $T$, the number of
partial sums \eqnok{thetaj} into which the second-stage problems are
divided.

\item {\bf max processors}. The number of workers requested.

\item {\bf average processors}. The average of the number of active
(nonsuspended) worker processors available for use by our problem
during the run.  Because of the dynamic nature of the Condor system,
the actual number of available processors fluctuates continually
during the run.

\item {\bf parallel efficiency}. The proportion of time for which
  worker processors were kept busy solving second-stage problems
  while they were owned by this run.

\item {\bf maximum number of cuts in the model}. The maximum number of
(partial) subgradients that are used to define the model function
during the course of the algorithm.

\item {\bf masterproblem solve time}. The total time spent solving the
master subproblem to generate new candidate iterates during the course of the
algorithm.

\item {\bf wall clock}. The total time (in minutes) between submission
of the job and termination.

\end{itemize}

\begin{table}
\vspace*{1.0in}
\centering
\begin{tabular}{|c|r|rrr|rrr|rr|r|}
\begin{rotate}{-45} run \end{rotate} & 
\begin{rotate}{-45} points evaluated \end{rotate} & 
\begin{rotate}{-45} $\sigma$  \end{rotate} & 
\begin{rotate}{-45} \# tasks ($C$) \end{rotate} & 
\begin{rotate}{-45} \# clusters ($T$) \end{rotate} & 
\begin{rotate}{-45} max. processors allowed \end{rotate} &
\begin{rotate}{-45} av. processors \end{rotate} & 
\begin{rotate}{-45} parallel efficiency \end{rotate} & 
\begin{rotate}{-45} max. \# cuts in model \end{rotate} & 
\begin{rotate}{-45} masterproblem solve time (min) \end{rotate} & 
\begin{rotate}{-45} wall clock time (min) \end{rotate} \\ \hline

ALS & 269 & $.5$ & 10 & 50 & 20 & 15 & %
.74 & 5491 & 26 & 368 \\
ALS & 275 & $.5$ & 25 & 50 & 50 & 21 & %
.90 & 5536 & 25 & 270 \\
ALS & 293 & $.5$ & 50 & 50 & 100 & 20 & %
.83 & 5639 & 27 & 329 \\
ALS & 270 & $.7$ & 10 & 50 & 20 & 12 & %
.79 & 5522 & 27 & 509 \\
ALS & 274 & $.7$ & 25 & 50 & 50 & 25 & %
.73 & 5550 & 25 & 281 \\
ALS & 282 & $.7$ & 50 & 50 & 100 & 26 & %
.81 & 5562 & 24 & 254 \\
ALS & 254 & $.85$ & 10 & 50 & 20 & 12 & %
.58 & 5496 & 22 & 575 \\
ALS & 276 & $.85$ & 25 & 50 & 50 & 19 & %
.57 & 5575 & 23 & 516 \\
ALS & 278 & $.85$ & 50 & 50 & 100 & 35 & %
.49 & 5498 & 25 & 260 \\
\hline

\end{tabular} 
\caption{SSN, with $N=10,000$ scenarios, Algorithm ALS.\label{tab.ssn.10k.exp2}}
\end{table}

Table~\ref{tab.ssn.10k.exp2} shows the results of a series of trials
of Algorithm ALS with three different values of $\sigma$ ($.5$, $.7$,
and $.85$) and three different choices for the number of chunks $C$
into which the second-stage solutions were divided (10, 25, and 50).
The number of clusters $T$ was fixed at 50, so that up to 50
cuts were generated at each iteration.  For $\sigma=.5$, the number of
values of $x$ for which second-stage evaluations are occurring at any
point in time ranged from 2 to 4 during the runs, while for
$\sigma=.85$, there were never more than 2 points being evaluated
simultaneously.

When these runs were performed, we were not able to obtain anything
approaching the requested number $2C$ of workers from the Condor pool.
As general trends, we see that the less synchronous variants (with
$\sigma = .5$ and $\sigma=.7$) tend to be faster than the more
synchronous variant (with $\sigma=.85$), except for the final run,
during which more processors were available.  Moreover, larger values
of $C$ also tend to produce faster runs.  We also note that the number
of iterations does not depend strongly on $\sigma$. We would not, of
course, expect $C$ to affect strongly the number of iterations, but
since it affects the manner in which the second-stage evaluation work
is distributed, we {\em would} expect it to affect the run time. Since
the number of workers available to us during this run was limited,
however, we did not see the full benefit of a finer-grained work
distribution ($C=50$), though the relatively low parallel efficiency
of the final run ($\sigma=.85$, $C=50$) indicates that the benefits of
more processors may not have been great in any case.

A note on typical task sizes: For $C=10$, a typical task required
about $50$-$280$ seconds on a typical worker machine available to us,
while for $C=50$, about $9$-$60$ seconds were required. The large
variation reflects the wide range in processing ability of the
machines available in a pool during a typical run. These numbers also
generally hold for the results in Tables~\ref{tab.ssn.10k.exp4.2} and
\ref{tab.ssn.10k.exp4.1}.

By comparing the results from Table~\ref{tab.ssn.10k.exp2} with those
reported in Tables~\ref{tab.ssn.10k.exp4.2} and
\ref{tab.ssn.10k.exp4.1}, we  verified that Algorithm
ALS was not as efficient on this problem as Algorithm TR and certain
variants of Algorithm ATR. One advantage, however, was that the
asymptotic convergence of ALS was quite fast. Having taken many
iterations to build up a model and return to a neighborhood of the
solution after having strayed far from it in early iterations, the
last three to four iterations home in rapidly from a relatively crude
approximate solution (a relative accuracy $(\cQ_{\rm min} -
m(x^{k+1})) / (1 + | \cQ_{\rm min}|)$ of between $.0006$ and $.0026$)
to a solution of high accuracy. 
\begin{table}
\vspace*{1.0in}
\centering
\begin{tabular}{|c|r|rrr|rrr|rr|r|}
\begin{rotate}{-45} run \end{rotate} & 
\begin{rotate}{-45} points evaluated \end{rotate} & 
\begin{rotate}{-45} $|\cB|$ ($K$) \end{rotate} & 
\begin{rotate}{-45} \# tasks ($C$) \end{rotate} & 
\begin{rotate}{-45} \# clusters ($T$) \end{rotate} & 
\begin{rotate}{-45} max. processors allowed \end{rotate} &
\begin{rotate}{-45} av. processors \end{rotate} & 
\begin{rotate}{-45} parallel efficiency \end{rotate} & 
\begin{rotate}{-45} max. \# cuts in model \end{rotate} & 
\begin{rotate}{-45} masterproblem solve time (min) \end{rotate} & 
\begin{rotate}{-45} wall clock time (min) \end{rotate} \\ \hline

TR & 48 & - & 10 & 100 & 20 & 19 & .21 & 4284 & 3 & 131 \\
TR & 72 & - & 10 & 50 & 20 & 19 & .26 & 3520 & 3 & 150  \\
TR & 39 & - & 25 & 100 & 25 & 22 & .49 & 3126 & 2 & 59 \\
TR & 75 & - & 25 & 50 & 25 & 23 & .48 & 3519 & 3 & 114  \\
TR & 43 & - & 50 & 100 & 50 & 42 & .52 & 3860 & 3 & 35  \\
TR & 61 & - & 50 & 50 & 50 & 44 & .53 & 3011 & 3 & 40  \\
\hline

ATR & 109 & 3 & 10 & 100 & 20 & 18 & .74 & 7680 & 9 & 107  \\
ATR & 121 & 3 & 10 & 50 & 20 & 19 & .66 & 4825 & 6 & 111  \\
ATR & 105 & 3 & 25 & 100 & 50 & 37 & .73 & 7367 & 8 & 49  \\
ATR & 113 & 3 & 25 & 50 & 50 & 41 & .60 & 4997 & 6 & 48  \\
ATR & 103 & 3 & 50 & 100 & 100 & 66 & .55 & 7032 & 9 & 29  \\
ATR & 129 & 3 & 50 & 50 & 100 & 66 & .59 & 5183 & 7 & 32  \\
\hline

ATR & 167 & 6 & 10 & 100 & 35 & 24 & .93 & 7848 & 13 & 99  \\
ATR & 209 & 6 & 10 & 50 & 35 & 22 & .89 & 5730 & 15 & 92  \\
ATR & 186 & 6 & 25 & 100 & 87 & 49 & .77 & 8220 & 14 & 53  \\
ATR & 172 & 6 & 25 & 50 & 87 & 49 & .80 & 5945 & 7 & 49 \\
ATR & 159 & 6 & 50 & 100 & 175 & 31 & .89 & 7092 & 11 & 65  \\
ATR & 213 & 6 & 50 & 50 & 175 & 40 & .88 & 6299 & 12 & 70  \\
\hline

ATR & 260 & 9 & 10 & 100 & 50 & 12 & .95 & 14431 & 35 & 267  \\
ATR & 286 & 9 & 10 & 50 & 50 & 23 & .90 & 6528 & 19 & 160  \\
ATR & 293 & 9 & 25 & 100 & 125 & 17 & .93 & 9911 & 30 & 232  \\
ATR & 377 & 9 & 25 & 50 & 125 & 15 & .96 & 7080 & 24 & 321  \\
ATR & 218 & 9 & 50 & 100 & 200 & 28 & .82 & 10075 & 25 & 101  \\
ATR & 356 & 9 & 50 & 50 & 200 & 23 & .93 & 6132 & 23 & 194  \\
\hline

ATR & 378 & 14 & 10 & 100 & 75 & 18 & .88 & 15213 & 77 & 302  \\
ATR & 683 & 14 & 10 & 50 & 75 & 14 & .98 & 8850 & 48 & 648  \\
ATR & 441 & 14 & 25 & 100 & 187 & 22 & .89 & 14597 & 61 & 312  \\
ATR & 480 & 14 & 25 & 50 & 187 & 20 & .94 & 8379 & 36 & 347  \\
ATR & 446 & 14 & 50 & 100 & 200 & 20 & .83 & 13956 & 64 & 331  \\
ATR & 498 & 14 & 50 & 50 & 200 & 22 & .94 & 7892 & 35 & 329   \\
\hline

\end{tabular} 
\caption{SSN, with $N=10,000$ scenarios, first trial, Algorithms TR and ATR.\label{tab.ssn.10k.exp4.2}}
\end{table}

\begin{table}
\vspace*{1.0in}
\centering
\begin{tabular}{|c|r|rrr|rrr|rr|r|}
\begin{rotate}{-45} run \end{rotate} & 
\begin{rotate}{-45} points evaluated \end{rotate} & 
\begin{rotate}{-45} $|\cB|$ ($K$) \end{rotate} & 
\begin{rotate}{-45} \# tasks ($C$) \end{rotate} & 
\begin{rotate}{-45} \# clusters ($T$) \end{rotate} & 
\begin{rotate}{-45} max. processors allowed \end{rotate} &
\begin{rotate}{-45} av. processors \end{rotate} & 
\begin{rotate}{-45} parallel efficiency \end{rotate} & 
\begin{rotate}{-45} max. \# cuts in model \end{rotate} & 
\begin{rotate}{-45} masterproblem solve time (min) \end{rotate} & 
\begin{rotate}{-45} wall clock time (min) \end{rotate} \\ \hline

TR & 47 & - & 10 & 100 & 20 & 17 & .24 & 3849 & 4 & 192  \\
TR & 67 & - & 10 & 50 & 20 & 13 & .34 & 3355 & 3 & 256 \\
TR & 47 & - & 25 & 100 & 25 & 18 & .49 & 3876 & 4 & 97 \\
TR & 57 & - & 25 & 50 & 25 & 18 & .40 & 2835 & 3 & 119 \\
TR & 42 & - & 50 & 100 & 50 & 30 & .22 & 3732 & 3 & 122 \\
TR & 65 & - & 50 & 50 & 50 & 31 & .25 & 3128 & 4 & 151 \\
\hline

ATR & 92 & 3 & 10 & 100 & 20 & 11 & .89 & 7828 & 9 & 125 \\
ATR & 98 & 3 & 10 & 50 & 20 & 11 & .84 & 4893 & 5 & 173 \\
ATR & 86 & 3 & 25 & 100 & 50 & 34 & .38 & 6145 & 5 & 70 \\
ATR & 95 & 3 & 25 & 50 & 50 & 32 & .41 & 4469 & 4 & 77 \\
ATR & 80 & 3 & 50 & 100 & 100 & 52 & .23 & 5411 & 5 & 80 \\
ATR & 131 & 3 & 50 & 50 & 100 & 59 & .47 & 4717 & 6 & 55 \\
\hline

ATR & 137 & 6 & 10 & 100 & 35 & 30 & .57 & 8338 & 12 & 84 \\
ATR & 200 & 6 & 10 & 50 & 35 & 26 & .60 & 5211 & 9 & 130 \\
ATR & 119 & 6 & 25 & 100 & 87 & 52 & .55 & 7181 & 7 & 44 \\
ATR & 199 & 6 & 25 & 50 & 87 & 58 & .48 & 5298 & 9 & 81 \\
ATR & 178 & 6 & 50 & 100 & 175 & 50 & .47 & 9776 & 15 & 77 \\
ATR & 240 & 6 & 50 & 50 & 175 & 61 & .64 & 5910 & 11 & 74 \\
\hline

ATR & 181 & 9 & 10 & 100 & 50 & 37 & .56 & 8737 & 15 & 96 \\
ATR & 289 & 9 & 10 & 50 & 50 & 19 & .93 & 7491 & 25 & 238 \\
ATR & 212 & 9 & 25 & 100 & 125 & 90 & .66 & 11017 & 21 & 45 \\
ATR & 272 & 9 & 25 & 50 & 125 & 65 & .45 & 6365 & 15 & 105 \\
ATR & 281 & 9 & 50 & 100 & 200 & 51 & .72 & 11216 & 34 & 88 \\
ATR & 299 & 9 & 50 & 50 & 200 & 26 & .83 & 7438 & 27 & 225 \\
\hline

ATR & 304 & 14 & 10 & 100 & 75 & 38 & .89 & 13608 & 43 & 129 \\
ATR & 432 & 14 & 10 & 50 & 75 & 42 & .95 & 7844 & 28 & 132 \\
ATR & 356 & 14 & 25 & 100 & 187 & 71 & .78 & 13332 & 48 & 111 \\
ATR & 444 & 14 & 25 & 50 & 187 & 45 & .89 & 7435 & 36 & 163 \\
ATR & 388 & 14 & 50 & 100 & 200 & 42 & .79 & 12302 & 52 & 192 \\
ATR & 626 & 14 & 50 & 50 & 200 & 48 & .81 & 7273 & 46 & 254  \\ 
\hline
\end{tabular} 
\caption{SSN, with $N=10,000$ scenarios, second trial, Algorithms TR and ATR.\label{tab.ssn.10k.exp4.1}}
\end{table}

We now turn to Tables~\ref{tab.ssn.10k.exp4.2} and
\ref{tab.ssn.10k.exp4.1}, which report on two sets of trials on the
same problem as in Table~\ref{tab.ssn.10k.exp2}. In these trials we
varied the following parameters:
\bi
\item {\bf basket size:} 
$K=1$ (synchronous TR) as well as $K=3,6,9,14$;

\item {\bf number of tasks:} 
$C=10,25,50$, as in Table~\ref{tab.ssn.10k.exp2};

\item {\bf number of clusters:} $T=50,100$.
\ei
The parameter $\sigma$ was fixed at $.7$ in all these runs.

The results in Table~\ref{tab.ssn.10k.exp4.2} were obtained with the
master processor running on an Intel Solaris machine, while
Table~\ref{tab.ssn.10k.exp4.1} was obtained with a Linux master.  In
both cases, the Condor pool that we tapped for worker processors was
identical. Therefore, it is possible to do a meaningful comparison
between each line of Table~\ref{tab.ssn.10k.exp4.1} and its
counterpart in Table~\ref{tab.ssn.10k.exp4.2}.  Conditions on the
Condor pool varied between and during each trial. This fact, combined
with the properties of the algorithm, resulted in large variability of
runtime from one trial to the next, as we discuss below.

The nondeterministic nature of the algorithms is evident in doing a
side-by-side comparison of the two tables. Even for synchronous TR,
the slightly different numerical values for function and subgradient
value returned by different workers in different runs results in
slight variations in the iteration sequence and therefore slight
differences in the number of iterations. For the asynchronous
Algorithm ATR, the nondeterminism is even more marked.  During the
basket-filling phase of the algorithm, computation of a new $x$ is
triggered when a certain proportion of tasks from a current value of
$x$ has been returned. On different runs, the tasks will be returned
in different orders, so the information used by the trust-region
subproblem \eqnok{trsub.atr1} in generating the new point will vary
from run to run, and the resulting iteration sequences will generally
show substantial differences.

The synchronous TR algorithm is clearly better than the ATR variants
with $K>1$ in terms of total computation, which is roughly
proportional to the number of iterations. In fact, the total amount of
work increases steadily with basket size.  Because of the decreased
synchronicity requirements and the greater parallelism obtained for
$K>1$, the wall clock times (last columns) do not follow quite the
same trend. The wall clock times for basket sizes $K=3$ and $K=6$ are
at least competitive with the results obtained for the synchronous TR
algorithm. The choice $K=6$ gave few of the fastest runs but did yield
consistent performance over all the different choices for the other
parameters, and under different Condor pool conditions.

The deleterious effects of synchronicity in Algorithm TR can be seen in
its poor performance on several instances, particularly during the
second trial. Let us compare, for instance, the entries in the two
tables for the variant of TR with $C=50$ and $T=100$. In the first
trial, this run used 42 worker processors on average and took 35
minutes, while in the second trial it used 30 workers on average and
required 122 minutes. The difference in runtime is too large to be
accounted for by the number of workers. Because this is a synchronous
algorithm, the time required for each iteration is determined by the
time required for the slowest worker to return the results of its
task. In the first trial, almost all tasks required between 6 and 35
seconds, except for a few iterations that contained tasks that took up
to 62 seconds. In the second trial, the slowest worker at each
iteration almost always required more than 60 seconds to complete its
task. We return to this point in discussing
Table~\ref{tab.ssn.10k.exp5} below.

Other general observations we can make are that 100 clusters give
almost uniformly better results in terms of wall clock time than 50
clusters, although the higher number results in a larger number of
cuts in the trust-region subproblems and an increased amount of time
on the master processor in solving these problems. The latter factor
is critical for $K=9$ and $K=14$, which do not compare
favorably with the smaller values of $K$ on this problem, even if many
more worker processors are available.  For the large basket sizes, the
loss of control induced by the increase in assynchronicity leads to a
significantly larger number of points that are evaluated.

In all cases, it takes some time for the model $m$ to become a good
enough approximation to $\cQ$ that it generates a step that meets the
trust-region acceptance criteria. The six TR runs in
Table~\ref{tab.ssn.10k.exp4.1}, for instance, required 18, 27, 16, 22,
16, and 26 trust-region subproblems to be solved, respectively, before
they stepped away from the initial point. (Note that, as expected, the
runs with $T=100$ required fewer such iterations than those with
$T=50$.) After the first step is taken, most steps are successful;
that is, the first minor iterate usually is accepted as the next major
iterate. Occasionally, two to four minor iterations are required
before the next major iteration is identified.  Similar behavior is
observed for the runs of ATR, except that successful iterations are
more widely spaced. For the first run with $K=6$ in
Table~\ref{tab.ssn.10k.exp4.1}, for instance, the $37$th solution of
\eqnok{trsub.atr1} yields the first successful step; then 36 of the
following 99 solutions of the subproblem yield successful steps.

\begin{table}
\vspace*{1.0in}
\centering
\begin{tabular}{|c|r|rrr|rrr|rr|r|}
\begin{rotate}{-45} run \end{rotate} & 
\begin{rotate}{-45} points evaluated \end{rotate} & 
\begin{rotate}{-45} $|\cB|$ ($K$) \end{rotate} & 
\begin{rotate}{-45} \# tasks ($C$) \end{rotate} & 
\begin{rotate}{-45} \# clusters ($T$) \end{rotate} & 
\begin{rotate}{-45} max. processors allowed \end{rotate} &
\begin{rotate}{-45} av. processors \end{rotate} & 
\begin{rotate}{-45} parallel efficiency \end{rotate} & 
\begin{rotate}{-45} max. \# cuts in model \end{rotate} & 
\begin{rotate}{-45} masterproblem solve time (min) \end{rotate} & 
\begin{rotate}{-45} wall clock time (min) \end{rotate} \\ \hline

TR & 47 & - & 25 & 100 & 25 & 23 & .49 & 4040 & 3 & 58 \\
TR & 44 & - & 25 & 100 & 25 & 21 & .31 & 3220 & 3 & 97 \\
TR & 45 & - & 25 & 100 & 25 & 20 & .23 & 3966 & 4 & 158 \\ \hline

TR & 51 & - & 50 & 100 & 50 & 37 & .33 & 4428 & 3 & 48 \\
TR & 51 & - & 50 & 100 & 50 & 45 & .14 & 4806 & 3 & 135 \\
TR & 46 & - & 50 & 100 & 50 & 41 & .15 & 3847 & 4 & 135 \\ \hline

ATR & 81 & 3 & 25 & 100 & 50 & 43 & .38 & 7451 & 6 & 64 \\
ATR & 81 & 3 & 25 & 100 & 50 & 39 & .41 & 6461 & 5 & 64 \\
ATR & 87 & 3 & 25 & 100 & 50 & 36 & .44 & 6055 & 8 & 66 \\ \hline

ATR & 106 & 3 & 50 & 100 & 100 & 84 & .28 & 8222 & 9 & 53 \\
ATR & 95  & 3 & 50 & 100 & 100 & 65 & .26 & 6786 & 7 & 64 \\
ATR & 94  & 3 & 50 & 100 & 100 & 23 & .44 & 6593 & 8 & 105 \\ \hline

ATR & 171 & 6 & 25 & 100 & 87 & 70 & .45 & 9173 & 19 & 61 \\
ATR & 135 & 6 & 25 & 100 & 87 & 61 & .39 & 7354 & 12 & 75 \\
ATR & 145 & 6 & 25 & 100 & 87 & 38 & .35 & 8919 & 16 & 146 \\ \hline

ATR & 177 & 6 & 50 & 100 & 175 & 87 & .41 & 9263 & 22 & 54 \\
ATR & 162 & 6 & 50 & 100 & 175 & 93 & .34 & 7832 & 18 & 66 \\
ATR & 159 & 6 & 50 & 100 & 175 & 39 & .27 & 8215 & 22 & 199 \\ \hline

\end{tabular} 
\caption{SSN final trial with best parameter combinations, $N=10,000$ scenarios, Algorithms TR and ATR.\label{tab.ssn.10k.exp5}}
\end{table}

In Table~\ref{tab.ssn.10k.exp5}, we took the most promising parameter
combinations from Tables~\ref{tab.ssn.10k.exp4.1} and
\ref{tab.ssn.10k.exp4.2} and ran three trials with each combination.
The Condor pool conditions varied widely during this trial, as can be
seen by the way that the average number of workers varies within each
group of three runs. For the asynchronous ATR runs, the differences in
wall clock times within each set of three runs usually can be
explained in terms of the varying number of workers available. (A
possible exception is the last line of the table, the third run of ATR
with $K=6$, $C=50$ and $T=100$, which took almost four times as long
as the first run while having only slightly fewer than half as many
processors. While the speed of machines available was roughly similar
between these runs, the third run was plagued with numerous
suspensions as the workers were reclaimed by their owners. Total time
that workers were suspended was over 23,000 seconds on the third run
and less than 2,800 seconds during the first run.)  On the other hand,
the variability in wall clock time between the six runs of the
synchronous TR algorithm was due not to the number of available
workers but rather to the synchronicity effect described above. In the
run reported in the first line of the table, for instance, the slowest
worker on any iteration typically took less than 65 seconds. In the
run reported on the third line, the time required by the slowest
worker varied significantly but was typically much longer, 150 seconds
and more.

\subsection{Larger Instances} \label{sec:results:large}

We also performed runs on several larger instances of SSN (with
$N=100,000$ scenarios) and on some very large instances
of the stormG2 problem, a cargo flight scheduling application described
by Mulvey and Ruszczy{\'n}ski~\cite{MulR95}.  Our interest
in this section is more in the sheer size of the problems that can be
solved using the algorithms developed for the computational grid
than with the relative performance of the algorithms with
different parameter settings.  

\begin{table}
\vspace*{1.0in}
\centering
\begin{tabular}{|c|r|rrr|rrr|rr|r|}
\begin{rotate}{-45} run \end{rotate} & 
\begin{rotate}{-45} points evaluated \end{rotate} & 
\begin{rotate}{-45} $|\cB|$ ($K$) \end{rotate} & 
\begin{rotate}{-45} \# tasks ($C$) \end{rotate} & 
\begin{rotate}{-45} \# clusters ($T$) \end{rotate} & 
\begin{rotate}{-45} max. processors allowed \end{rotate} &
\begin{rotate}{-45} av. processors \end{rotate} & 
\begin{rotate}{-45} parallel efficiency \end{rotate} & 
\begin{rotate}{-45} max. \# cuts in model \end{rotate} & 
\begin{rotate}{-45} masterproblem solve time (min) \end{rotate} & 
\begin{rotate}{-45} wall clock time (min) \end{rotate} \\ \hline
ATR & 177 & 3 & 100 & 100 & 200 & 38 & .52 & 10558 & 47 & 1357 \\ 
\hline
\end{tabular}
\caption{SSN, with $N=100,000$ scenarios.\label{tab.ssn.100k}}
\end{table}

Table~\ref{tab.ssn.100k} shows results for a sampled instance of SSN
with $N=100,000$ scenarios, which is a linear program with
approximately $1.75 \times 10^7$ constraints and $7.06 \times 10^7$
variables. This run was performed at a time when not many machines
were available, and many suspensions occurred during the run. We chose
$T=100$ chunks per evaluation and found that most tasks required
between 41 and 300 seconds on the workers, with a few task times of
more than 500 seconds. (The benchmarks indicated that the worker speed
varied over a factor of 7.)  A total of 77 different workers were used
during the run, though the average number of nonsuspended workers
available at any time was only 39. In fact, at any given point in the
computation there were an average of 7 workers assigned to this task
that were suspended. Still, a result was obtained in about 22 hours.

\begin{table}
\vspace*{1.0in}
\centering
\begin{tabular}{|c|r|rrr|rrr|rr|r|}
\begin{rotate}{-45} run \end{rotate} & 
\begin{rotate}{-45} points evaluated \end{rotate} & 
\begin{rotate}{-45} $|\cB|$ ($K$) \end{rotate} & 
\begin{rotate}{-45} \# tasks ($C$) \end{rotate} & 
\begin{rotate}{-45} \# clusters ($T$) \end{rotate} & 
\begin{rotate}{-45} max. processors allowed \end{rotate} &
\begin{rotate}{-45} av. processors \end{rotate} & 
\begin{rotate}{-45} parallel efficiency \end{rotate} & 
\begin{rotate}{-45} max. \# cuts in model \end{rotate} & 
\begin{rotate}{-45} masterproblem solve time (min) \end{rotate} & 
\begin{rotate}{-45} wall clock time (min) \end{rotate} \\ \hline
TR  & 17 & -   & 125 & 125  & 250 & 106 & .55 & 2310 & 0.5 & 146  \\ %
ATR & 25 & 3  & 125 & 125 & 250 & 106 & .90 & 3292 & 0.5 & 116 \\ \hline %
\end{tabular} 
\caption{stormG2, with $N=250000$ scenarios. \label{tab.storm.250k}}
\end{table}

In the stormG2 problem of Mulvey and Ruszczy{\'n}ski~\cite{MulR95}, the
first-stage problem contained 121 variables, while each second-stage
problem contained 1259 variables.  We considered first a sampled
approximation of this problem with 250000 scenarios, which resulted
in a linear program with $1.32 \times 10^8$ constraints and $315 \times 10^8$
unknowns.  Results are shown in Table~\ref{tab.storm.250k}. The
algorithm was started at a solution of a sampled instance with fewer
scenarios and was quite close to optimal. The objective function at
the initial point was approximately $15499595.1$, compared with an
optimal value of $15499591.9$ achieved by Algorithm TR. In fact, the
TR algorithm takes only one major iteration---it accepts the 16th
minor iteration as the first major iterate $x^1$. The ATR variant does
not take even one step---it terminates after determining that the
initial point $x^0$ is optimal to within the given convergence
tolerance.  Although we requested 250 processors, an average of only
106 were available during the time that we performed these two test
runs. The second run is able to utilize these to high efficiency, as
the second-stage workload can be divided into a large number of chunks
and very little time is spent in solving the trust-region subproblem.

\begin{table}
\vspace*{1.0in}
\centering
\begin{tabular}{|c|r|rrr|rrr|rr|r|}
\begin{rotate}{-45} run \end{rotate} & 
\begin{rotate}{-45} points evaluated \end{rotate} & 
\begin{rotate}{-45} $|\cB|$ ($K$) \end{rotate} & 
\begin{rotate}{-45} \# tasks ($C$) \end{rotate} & 
\begin{rotate}{-45} \# clusters ($T$) \end{rotate} & 
\begin{rotate}{-45} max. processors allowed \end{rotate} &
\begin{rotate}{-45} av. processors \end{rotate} & 
\begin{rotate}{-45} parallel efficiency \end{rotate} & 
\begin{rotate}{-45} max. \# cuts in model \end{rotate} & 
\begin{rotate}{-45} masterproblem solve time (hr) \end{rotate} & 
\begin{rotate}{-45} wall clock time (hr) \end{rotate} \\ \hline
ATR & 28 & 4 & 1024 & 1024 & 800 & 433 & .668 & 39647 & 1.9 & 31.9 \\ \hline
\end{tabular}
\caption{stormG2, with $N=10^7$ scenarios.\label{tab.storm.1e7}}
\end{table}

Finally, we report on a very large sampled instance of stormG2 with
$N=10^7$ scenarios, an instance whose deterministic equivalent is a
linear program with $9.85 \times 10^8$ constraints and $1.26 \times
10^{10}$ variables.  Performance is profiled in
Table~\ref{tab.storm.1e7}.

We used the tighter convergence tolerance $\epstol = 10^{-6}$ for this
run. The algorithm took successful steps at iterations 28, 34, 37, and
38, the last of these being the final iteration. The first evaluated
point had a function value of 
$15526740$, compared with a value of
$15498842$ at the final iteration. 

For this run, we augmented the Wisconsin Computer Science Condor pool with
machines from Georgia Tech, the University of New Mexico, the Italian
National Institute of Physics (INFN), the NCSA at the University of Illinois,
and the IEOR Department at Columbia, the Albu, and the Wisconsin
engineering Department.  Table~\ref{bigstorm.tab} shows
the number and type of processors available at each of these
locations. 
In contrast to the other runs
reported here, we used the ``MW-files'' implementation of MW, the
variant that uses shared files to perform communication between master
and workers rather than Condor-PVM.

\begin{table}
\centering
\begin{tabular}{|c|c|c|} \hline
Number & Type & Location \\ \hline
184 & Intel/Linux & Argonne \\ \hline
254  & Intel/Linux & New Mexico \\ \hline
36  & Intel/Linux & NCSA \\ \hline
265 & Intel/Linux & Wisconsin \\
88 & Intel/Solaris & Wisconsin \\
239 & Sun/Solaris & Wisconsin \\ \hline
124 & Intel/Linux & Georgia Tech  \\
90  & Intel/Solaris & Georgia Tech  \\
13 & Sun/Solaris & Georgia Tech \\ \hline
9   & Intel/Linux & Columbia U.  \\
10  & Sun/Solaris & Columbia U.  \\ \hline
 33  & Intel/Linux & Italy (INFN)  \\ \hline \hline
1345 & & \\ \hline
\end{tabular}
\caption{Machines available for stormG2, with $N=10^7$
scenarios.\label{bigstorm.tab}}
\end{table}

\begin{figure}
\centering
\epsfig{figure=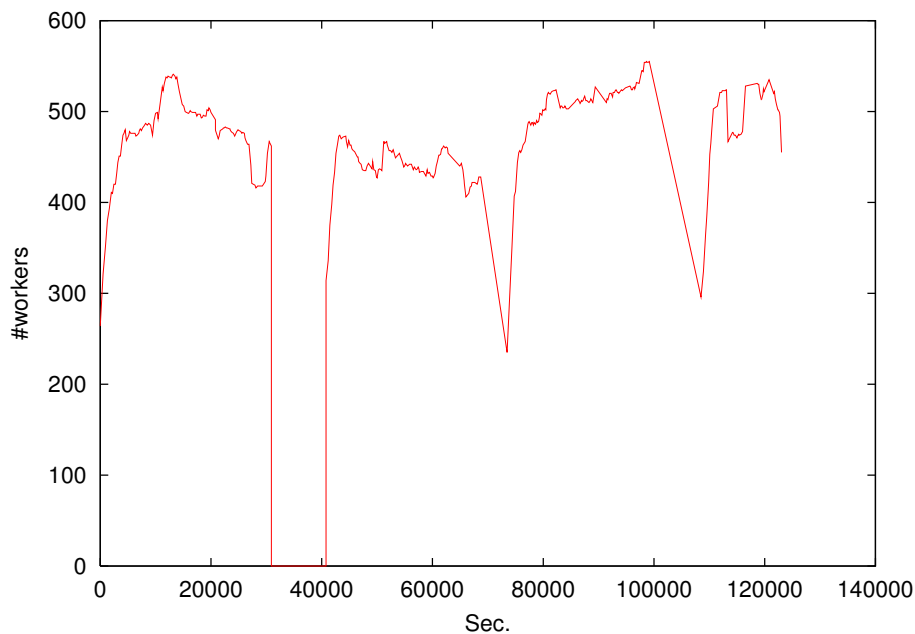,angle=270,width=\linewidth}
\caption{Number of workers used for stormG2, with $N=10^7$ scenarios.\label{bigstorm-workers.fig}}
\end{figure}

The job ran for a total of almost 32 hours.  The number of workers
being used during the course of the run is shown in
Figure~\ref{bigstorm-workers.fig}.  The job was stopped after
approximately 8 hours and was restarted manually from a checkpoint
about 2 hours later.  It then ran for approximately 24 hours to
completion.  The number of workers dopped off significantly on two
occasions.  The drops were due to the master processor ``blocking'' to
solve a difficult master problem and to checkpoint the state of the
computation.  During this time the worker processors were idle, and
MW decided to release a number of the processors rather than have them
sit idle.

As noted in Table~\ref{tab.storm.1e7}, an average of 433 workers were
present at any given point in the run. The computation used a maximum
of 556 workers, and there was a ratio of 12 in the speed of the
slowest and fastest machines, as determined by the benchmarks. A total
of 40837 tasks were generated during the run, representing $3.99
\times 10^8$ second-stage linear programs. (At this rate, an average
of 3472 second-stage linear programs were being solved per second
during the run.) The average time to solve a task was 774 seconds.
The total cumulative CPU time spent by the worker pool was 9014 hours,
or just over one year of computation.

\section{Conclusions}

We have described L-shaped and trust-region algorithms for solving the
two-stage stochastic linear programming problem with recourse, and
derived asynchronous variants suitable for parallel implementation on
distributed heterogeneous computational grids. We prove convergence
results for the trust-region algorithms. Implementations based on the
MW library and the Condor system are described, and we report on
computational studies using different algorithmic parameters under
different pool conditions.  Becasue of the dynamic nature of the
computational pool, it is impossible to arrive at a ``best''
configuration or set of algorithmic parameters for all instances.
Instead, it may be important to adjust the algorithm parameters
dynamically; we suggest this as a line of future research.  Finally,
we report on the solution of some large sampled instances of problems
from the literature, including an instance of the stormG2 problem
whose deterministic equivalent has more than $10^{10}$ unknowns.
Since the use of the computational grid has the greatest benefit on
problems that require large amounts of computation, the algorithms
developed here are best suited to larger (multistage) problems or
incorporated into a sample average approximation approach (see Shapiro and Homem-de-Mello~\cite{ShaH01}.

\section*{Acknowledgments}

This research was supported by the Mathematics, Information, and
Computational Sciences Division subprogram of the Office of Advanced
Scientific Computing Research, U.S. Department of Energy, under
Contract W-31-109-Eng-38.  We also acknowledge the support of the
National Science Foundation, under Grant CDA-9726385.  We would also
like to acknowledge the IHPCL at Georgia Tech, which is supported by a
grant from Intel; the National Computational Science Alliance under
grant number MCA00N015N for providing resources at the University of
Wisconsin, the NCSA SGI/CRAY Origin2000, and the University of New
Mexico/Albuquerque High Performance Computing Center AltaCluster; and
the Italian Istituto Nazionale di Fisica Nucleare (INFN) and Columbia
University for allowing us access to their Condor pools.

We are grateful to Alexander Shapiro and Sven Leyffer for discussions
about the algorithms presented here.

\bibliographystyle{plain}
\bibliography{refs}

\begin{thebibliography}{10}

\bibitem{BahDGV95}
O.~Bahn, O.~{du Merle}, {J.-L.} Goffin, and J.~P. Vial.
\newblock A cutting-plane method from analytic centers for stochastic
  programming.
\newblock {\em Mathematical Programming, Series B}, 69:45--73, 1995.

\bibitem{Ben62}
J.~F. Benders.
\newblock Partitioning procedures for solving mixed variable programming
  problems.
\newblock {\em Numerische Mathematik}, 4:238--252, 1962.

\bibitem{BirDGGKW87}
J.~R. Birge, M.~A.~H. Dempster, H.~I. Gassmann, E.~A. Gunn, and A.~J. King.
\newblock A standard input format for multiperiod stochastic linear programs.
\newblock {\em COAL Newsletter}, 17:1--19, 1987.

\bibitem{BirDHS98}
J.~R. Birge, C.~J. Donohue, D.~F. Holmes, and O.~G. Svintsiski.
\newblock A parallel implementation of the nested decomposition algorithm for
  multistage stochastic linear programs.
\newblock {\em Mathematical Programming}, 75:327--352, 1996.

\bibitem{BirL97}
J.~R. Birge and R.~Louveaux.
\newblock {\em Introduction to Stochastic Programming}.
\newblock Springer, New York, 1997.

\bibitem{BirQ88}
J.~R. Birge and L.~Qi.
\newblock Computing block-angular {Karmarkar} projections with applications to
  stochastic programming.
\newblock {\em Management Science}, 34:1472--1479, 1988.

\bibitem{BurF93}
J.~V. Burke and M.~C. Ferris.
\newblock Weak sharp minima in mathematical programming.
\newblock {\em SIAM Journal on Control and Optimization}, 31:1340--1359, 1993.

\bibitem{FraGV00}
E.~Frangi{\`e}re, J.~Gondzio, and J.-P. Vial.
\newblock Building and solving large-scale stochastic programs on an affordable
  distributed computing system.
\newblock {\em Annals of Operations Research}, 2000.
\newblock To appear.

\bibitem{GasS97}
H.~I. Gassmann and E.~Schweitzer.
\newblock A comprehensive input format for stochastic linear programs.
\newblock Working Paper WP-96-1, School of Business Administration, Dalhousie
  University, Halifax, Canada, December 1997.

\bibitem{PVMbook}
A.~Geist, A.~Beguelin, J.~Dongarra, W.~Jiang, R.~Manchek, and V.~Sunderam.
\newblock {\em {PVM}: Parallel Virtual Machine}.
\newblock The {MIT} Press, Cambridge, MA, 1994.

\bibitem{GonV00}
J.~Gondzio and J.-P Vial.
\newblock Warm start and $\epsilon$-subgradients in the cutting plane scheme
  for block-angular linear programs.
\newblock {\em Computational Optimization and Applications}, 14:17--36, 1999.

\bibitem{GouKLY00}
J.-P. Goux, S.~Kulkarni, J.~T. Linderoth, and M.~E. Yoder.
\newblock An enabling framework for master-worker applications on the
  computational grid.
\newblock In {\em Proceedings of the Ninth {IEEE} Symposium on High Performance
  Distributed Computing}, 2000.

\bibitem{GouLY00}
J.-P. Goux, J.~T. Linderoth, and M.~E. Yoder.
\newblock Metacomputing and the master-worker paradigm.
\newblock Preprint ANL/MCS-P792-0200, Mathematics and Computer Science
  Division, Argonne National Laboratory, 2000.

\bibitem{HirL93}
{J.-B.} {Hiriart-Urruty} and C.~{Lemar{\'e}chal}.
\newblock {\em Convex Analysis and Minimization Algorithms {II}}.
\newblock Comprehensive Studies in Mathematics. Springer-Verlag, 1993.

\bibitem{Hol97}
1997.
\newblock {\tt http://www-personal.umich.edu/\~jrbirge/dholmes/SPTSlists.html}.

\bibitem{Kiw90}
K.~C. Kiwiel.
\newblock Proximity control in bundle methods for convex nondifferentiable
  minimization.
\newblock {\em Mathematical Programming}, 46:105--122, 1990.

\bibitem{condor}
M.~Livny, J.~Basney, R.~Raman, and T.~Tannenbaum.
\newblock Mechanisms for high throughput computing.
\newblock {\em SPEEDUP}, 11, 1997.
\newblock Available from {\tt http://www.cs.wisc.edu/condor/doc/htc\_mech.ps}.

\bibitem{Man69}
O.~L. Mangasarian.
\newblock {\em Nonlinear Programming}.
\newblock McGraw-Hill, New York, 1969.

\bibitem{MulR95}
J.~M. Mulvey and A.~Ruszczy{\'n}ski.
\newblock A new scenario decomposition method for large scale stochastic
  optimization.
\newblock {\em Operations Research}, 43:477--490, 1995.

\bibitem{Roc70}
R.~T. Rockafellar.
\newblock {\em Convex Analysis}.
\newblock Princeton University Press, Princeton, 1970.

\bibitem{Rus86}
A.~Ruszczy{\'n}ski.
\newblock A regularized decomposition for minimizing a sum of polyhedral
  functions.
\newblock {\em Mathematical Programming}, 35:309--333, 1986.

\bibitem{Rus93}
A.~Ruszczy{\'n}ski.
\newblock Parallel decomposition of multistage stochastic programming problems.
\newblock {\em Mathematical Programming}, 58:201--228, 1993.

\bibitem{SenDC94}
S.~Sen, R.~D. Doverspike, and S.~Cosares.
\newblock Network planning with random demand.
\newblock {\em Telecommunications Systems}, 3:11--30, 1994.

\bibitem{ShaH01}
Alexander Shapiro and Tito {Homem-de-Mello}.
\newblock On the rate of convergence of optimal solutions of {Monte Carlo}
  approximations of stochastic programs.
\newblock {\em SIAM Journal on Optimization}, 11(1):70--86, 2001.

\bibitem{VanW69}
R.~{Van Slyke} and {R.J-B.} Wets.
\newblock L-shaped linear programs with applications to control and stochastic
  programming.
\newblock {\em {SIAM} Journal on Applied Mathematics}, 17:638--663, 1969.

\bibitem{soplex}
R.~Wunderling.
\newblock {\em Paralleler und Objektorientierter Simplex-Algorithmus}.
\newblock PhD thesis, Konrad-Zuse-Zentrum f{\"u}r Informationstechnik, Berlin,
  1996.

\end{thebibliography}

\end{document}